\documentclass[envcountsame,envcountsect]{svmult}

\usepackage{makeidx}         % allows index generation
\usepackage{graphicx}        % standard LaTeX graphics tool
\usepackage{multicol}        % used for the two-column index
\usepackage[bottom]{footmisc}% places footnotes at page bottom
\usepackage{amssymb}
\def \sl {\mathfrak{sl}}

\begin{document}

%%%%%%%%%%%%%%%%%%%%%%    Title    %%%%%%%%%%%%%%%%%%%%%
\title*{Poisson and symplectic functions\\ in Lie algebroid theory}
\titlerunning{Poisson and symplectic functions}
\author{Yvette Kosmann-Schwarzbach\inst{}}
\institute{Centre de Math\'ematiques Laurent Schwartz\\\'Ecole
  Polytechnique\\ 91128 Palaiseau, France\\
\texttt{yks@math.polytechnique.fr}}
%%%%%%%%%%%%%%%%%%%%%%

\maketitle

\hfill {\it{For Murray Gerstenhaber and Jim
    Stasheff, in admiration and respect}}

\begin{abstract}
Emphasizing the role of Gerstenhaber algebras
and of higher derived brackets in the theory of Lie algebroids, 
we show that the several Lie algebroid brackets 
which have been introduced in the recent literature can all be 
defined in terms of Poisson and pre-symplectic 
functions in the sense of Roytenberg and Terashima.
We prove that in
this very general framework there exists a 
one-to-one correspondence between non-degenerate Poisson
functions and symplectic functions. We also determine the 
differential associated to a Lie algebroid structure obtained by
twisting a structure with background
by both a Lie bialgebra action and a Poisson bivector. 
\end{abstract}

\section*{Introduction}

Towards 1958, Ehresmann \cite{E1} introduced the idea of differentiable
categories, of which 
the differentiable groupoids, now called 
Lie groupoids, are an example, and he developed this theory further 
in the 1960's  
\cite{E2}. At the end of the decade,
Pradines introduced the corresponding infinitesimal objects which he called
Lie algebroids \cite{P}.
The theory of Lie algebroids, which has since been developed by many
authors, and in particular by Mackenzie \cite{M1} \cite{M2},
encompasses both differential
geometry -- because the tangent bundle of a smooth manifold is the 
prototypical Lie algebroid --, and Lie algebra theory -- because the Lie
algebras are Lie algebroids whose base
manifold is a singleton --, while other examples of Lie algebroids 
occur in the theory of foliations 
(see, {\it{e.g.}}, \cite{MM}) and in Poisson geometry
\cite{CDW} \cite{M2}. The 
corresponding, purely algebraic concept, called   
pseudo-Lie algebras \cite{M0} or Lie-Rinehart algebras \cite{H1},
among many other names, dates back to Jacobson \cite{J}, as has been 
observed in \cite{KMc}.

While the structure of what is now called a {\it {Gerstenhaber algebra}}
appeared in the work of Murray Gerstenhaber on the Hochschild
cohomology of associative algebras \cite{G}, it became clear  
in the work of Koszul \cite{Ko} and of 
many other authors \cite{KM} \cite{X}
that Gerstenhaber algebras play an
essential role in the theory of Lie algebroids. Whenever a vector
bundle has the structure of a Lie algebroid, the linear space of 
sections of its exterior algebra bundle is a Gerstenhaber algebra, the
prototypical example of which is the linear space of fields of multivectors
equipped with the Schouten--Nijenhuis bracket on any smooth manifold.
The close relationship between Poisson geometry and Lie algebroid
theory appears clearly in the concept of 
a Lie bialgebroid defined by Mackenzie and Xu \cite{MX} 
as the infinitesimal object of 
a Poisson groupoid, and characterized in terms of derivations in
\cite{yks95}. For any Poisson manifold $M$ with tangent bundle $TM$, the pair
$(TM,T^*M)$ is a Lie bialgebroid, while the Lie
bialgebroids over a point are Drinfeld's Lie bialgebras 
of Poisson-Lie group theory \cite{D}.

When passing from the case of Lie bialgebras to
that of the  
Lie-quasi\break bialgebras \cite{D2}, or their dual version, the 
quasi-Lie bialgebras, or 
the more general case of proto-bialgebras  \cite{yks92}\footnote{In
  \cite{yks92} \cite{BK}, Lie-quasi bialgebras were
  called Jacobian quasi-bialgebras, and quasi-Lie bialgebras were
  called co-Jacobian quasi-bialgebras. We also point out that,
 in the translation of Drinfeld's original
paper \cite{D2}, the term ``quasi-Lie bialgebra'' is used for 
what we call Lie-quasi bialgebra.
Proto-bialgebras were
  introduced in \cite{yks92} where they were called proto-Lie-bialgebras, to
  distinguish them from the associative version of this notion.},
 {\it {higher structures}},
in the sense of Jim Stasheff \cite{Sta}, appear.
The associated algebra is not
a Gerstenhaber algebra but only a Gerstenhaber algebra up to homotopy,
but with all $n$-ary brackets beyond the third vanishing (see
\cite{H2} \cite{B1} \cite{B2}).
The analogous theory generalizing Lie algebroids 
was developed by Roytenberg \cite{R} and, 
more recently, by Terashima
\cite{T}. Their articles form the basis of the present 
exposition\footnote{There are
some changes in the notations. 
In particular the notations $\phi$ and $\psi$ used by Roytenberg in
  \cite{R} are exchanged in order to return to the conventions of
  \cite{yks92} \cite{BK} \cite{yks05}.}. 

\medskip

The concept of {\it twisting} for proto-bialgebroids 
was defined by Roytenberg \cite{R} as a generalization of
  the twisting  
of proto-bialgebras introduced in \cite{yks92}, itself a generalization of the
twisting of Lie bialgebras defined by Drinfeld in the
theory of the semi-classical limit of the quasi-Hopf algebras \cite{D2}, 
while the concept of 
{\it Poisson function}, which was already implicit in \cite{R}, has
now been formally
introduced by Terashima
in \cite{T}, with interesting applications which we review and
  develop here.
Poisson functions generalize both Poisson structures on
manifolds and triangular $r$-matrices on Lie algebras, and, more generally,
Poisson structures on Lie algebroids
as well as their twisted versions  
(see \cite{LWX} \cite{R} \cite{T}). 

\medskip

The cohomological approach to Lie algebroid theory arose from the
viewpoint developed for Lie bialgebras by 
Lecomte and Roger \cite{LR}, itself based on the even Poisson bracket
introduced by Kostant and Sternberg in \cite{KSt}\footnote{Even 
  Poisson brackets had already 
appeared in the context of the quantization of systems
  with constraints in the work of Batalin,
  Fradkin and Vilkovisky. See \cite{Stasheff} and references therein.}. 
In \cite{yks92}, we extended this approach to the 
Lie-quasi bialgebras defined by Drinfeld \cite{D2}, and we introduced
the dual objects and the more general notion of proto-Lie bialgebra,
encompassing both the Lie-quasi bialgebras and their 
duals. In \cite{R} Roytenberg 
extended the cohomological approach to Lie bialgebra theory to the
``oid'' case 
by combining the supermanifold approach due to 
Vaintrob \cite{V} and  
T. Voronov (see \cite{Vo} citing earlier publications)
with the results of \cite{yks92}.

The preprint that 
Terashima communicated to me in 2006 \cite{T} goes further along the same
lines and provides a beautiful unification
of results in both recent \cite{BC} \cite{BCS} and not so recent papers
\cite{L}, 
showing that they
are special cases of a general construction of Lie algebroid
structures obtained by twisting certain
basic structures.

\medskip

The main features of this paper are the following.
Our first Section deals with the general definition of
a {\it{structure}} on a vector bundle, $V$. 
The basic tool for the study of the properties of
``structures'' is the {\it{big
  bracket}}, denoted by $\{~,~\}$, the bigraded even Poisson bracket 
which is the canonical
Poisson bracket on the cotangent bundle
of the supermanifold $\varPi V$, {\it{i.e.}}, 
$V$ with reversed parity on the
fibers,
which, on vector-valued
forms or $1$-form-valued multivectors, coincides with the
Nijenhuis-Richardson bracket up to sign. The ``structures'' are cubic functions
on this cotangent bundle whose Poisson square vanishes.   
Vector bundles equipped with a ``structure'' generalize the 
Lie, Lie-quasi and quasi-Lie
bialgebroids, in particular the Lie bialgebras.

In Section \ref{2}, we introduce the dual notions of {\it{twisting by a
    bivector}}
and {\it{twisting by a $2$-form}}, and we define the {\it{Poisson
    functions}} and the {\it{pre-symplectic functions}} with respect
to a given structure. 
Such bivectors (resp., $2$-forms) give rise by twisting to quasi-Lie
(resp., Lie-quasi) bialgebroids.
We show that the {\it{twist of Lie-quasi bialgebras}} 
in the sense of Drinfeld \cite{D2} and the {\it{twisted 
Poisson structures}} on manifolds, introduced by Klim\v c\'{\i}k and
Strobl in  \cite{KS} (under the name WZW-Poisson structures) and
studied by \v Severa and Weinstein in \cite{SW} (where they are called
Poisson structures with background),
are both particular cases of the general notion of a twisted structure. 

In Section \ref{3} we prove that the graphs of Poisson functions and
of pre-symplectic functions are
Dirac sub-bundles of the 
Courant algebroid $V \oplus V^*$, which is the ``double'' of $V$. 

The aim of Section \ref{4} is to prove Theorem \ref{symplectic}, which 
states that non-degenerate Poisson
functions are in one-to-one correspondence with symplectic functions,
a generalization of the well-known fact that a non-degnerate
bivector on a manifold defines a Poisson structure if and only if its
inverse is a closed $2$-form.
We believe that this theorem had not yet been proved in so general a form.

In Section \ref{another}, we study the case where 
a Poisson function involves
both a Poisson structure on a manifold $M$
in the ordinary sense and a Lie algebra action on this manifold. 
In the general case, with non trival Lie-quasi\break bialgebra
actions and 
background $3$-forms on the manifold, we determine explicit expressions
for the bracket and the differential thus defined. In fact,
the twisting of a structure on a vector bundle $V$ 
by a Poisson function gives rise to  a Lie algebroid structure on the
dual vector bundle
$V^*$ and, dually, to a differential on the sections of
$\wedge^\bullet V$, the exterior algebra bundle of $V$.
In particular cases, we recover the
brackets on vector bundles of the form $T^*M \times {\mathfrak g}$
which were associated to Poisson actions of Poisson-Lie groups
on Poisson manifolds by Lu in \cite{L} and, more generally, to 
quasi-Poisson
$G$-manifolds in the sense\break of \cite{AKM} by Bursztyn
and Crainic in \cite{BC}, and to
quasi-Poisson $G$-spaces in the sense of \cite{AK} by Bursztyn, Crainic
and \v Severa in \cite{BCS}.
This approach gives an immediate proof that these brackets satisfy the
Jacobi identity and are indeed Lie algebroid brackets.
The formulas for the differential in the general case are, to the best
of our knowledge, new.

\section{Definition of structures}\label{1} 
\subsection{Towards a unification}
It was already clear in the theory of Lie bialgebras that the 
``big bracket'' was the appropriate tool for their study. 
Roytenberg 
extended the definition and the use of the big bracket to the case of
Lie algebroids \cite{R}, and
Terashima's article \cite{T} proves additional results, by
suitably twisting certain
basic structures.

\subsection{The big bracket}
Consider the
bigraded supermanifold $X = T^*\varPi V$, 
where $V$ is a vector bundle over a manifold $M$, and where $\varPi$ denotes
the change of parity of the fibers. Then $X$ is 
canonically equipped with an even Poisson bracket \cite{KSt}, 
the Poisson structure on $X$ actually being symplectic.
This Poisson bracket, called the \emph{big bracket},
is here denoted by $\{~,~\}$. 
The algebra $\mathcal F$  of smooth functions 
on $X$ is bigraded in the following way. 
If $(x^{i},\xi^{a})$ are 
local coordinates on $\varPi V$ ($i= 1, \ldots, \dim M$, \, $a= 1,
\ldots, {\mathrm{rank}} \, V$), we denote by  
$(x^{i},\xi^{a}, p_i, \theta_a)$ the corresponding local coordinates
on $T^{*}\varPi V$, and we assign them the bidegrees $(0,0)$,  $(0,1)$,
$(1,1)$ and $(1,0)$, respectively.
An element of $\mathcal F$ of bidegree $(k,\ell)$, 
with $k \geq 0$ and $\ell \geq 0$, 
is said to be of {\it{shifted bidegree}} $(p,q)$
when $p= k-1$ and $q= \ell-1$ ($p \geq -1$ and $q \geq -1$), whence
the table
$$
\begin{array}{ccccc}
 x^{i} & \xi^{a}& p_i & \theta_a & \\ 
(0,0) & (0,1) & (1,1) & (1,0) & \mathrm{bidegree}\\
(-1,-1) & (-1,0)& (0,0) &  (0,-1) & 
\quad \mathrm{shifted} \, \, \mathrm{bidegree}\\ 
\end{array}
$$
The total
degree (resp., total shifted degree) will be called, for short, the degree
(resp., shifted degree).
The big bracket  is of shifted bidegree $(0,0)$, and it
satisfies
$$
\{ x^{i},p_{j} \} = \delta^{i}_{j} = - \{ p_{j}, x^{i} \} \ , \quad
\quad \{ \xi^{a}, \theta_{b}\} =
\delta^{a}_{b} = \{\theta_{b}, \xi^{a} \} \ .
$$

\subsection{Definition of structures}
As in \cite{V} \cite{R} \cite{Vo} (also see
\cite{yks05}) we consider functions on $X$ that define bialgebroid
structures or generalizations thereof on $(V,V^*)$.
See \cite{yks92} \cite{BK} \cite{R}
for proofs of the statements in this section.
\begin{definition} 
\label{structuredef}
A {\emph {structure}} on $V$ is a 
homological function on $X$ of 
degree~$3$, {\it {i.e.}}, an element $S \in \mathcal F$ of
shifted degree $1$ such
that $\{S,S\} = 0$.
\end{definition}

Let 
\begin{equation}\label{structure}
S= \phi + \gamma + \mu + \psi \ 
\end{equation} 
in the notations of
 \cite{yks92} and \cite{BK}. Then,

\noindent $\bullet$ $\phi$, of shifted bidegree $(2,-1)$, is a $3$-form on $V^*$, 
$$\phi = \frac{1}{6} \phi ^{abc}
 \theta_a \theta _b \theta _c \ ,
$$
 $\bullet$ $\gamma$, of shifted bidegree 
$(1,0)$, defines an anchor, $a^* : V^* \to TM$, and a bracket on $V^*$,
$$ 
\gamma = (a^*)^{ib} p_i \theta_b + \frac{1}{2} \gamma ^{bc}_a \theta _b
\theta_c \xi^a \ ,
$$
 $\bullet$ $\mu$, of shifted bidegree $(0,1)$, defines an anchor, $a_*
 : V \to TM$, and a bracket on $V$,
$$ 
\mu = (a_{*})_b ^i p_i \xi^b + \frac{1}{2} \mu ^{a}_{bc} \theta _a
\xi^b \xi^c \ ,
$$
 $\bullet$ $\psi$, of shifted bidegree $(-1,2)$, is a $3$-form on $V$,
$$\psi = \frac{1}{6} \psi_ {abc}
 \xi^a \xi^b \xi^c \ .
$$
Then
$S$ is a structure if and only if
$$
\left\lbrace
\begin{array}{ll}
\vspace{.1cm}
&\frac{1}{2} \{ \mu,\mu\} + \{ \gamma ,\psi\} = 0 \ ,\\
\vspace{.1cm}
& \{\mu ,\gamma\} + \{ \phi ,\psi\} =0 \ ,\\
\vspace{.1cm}
& \frac{1}{2} \{\gamma ,\gamma\} + \{ \mu ,\phi\} =0 \ ,\\
\vspace{.1cm}
&\{ \mu ,\psi\} = 0 \ ,\\
\vspace{.1cm}
&\{ \gamma , \phi\} = 0 \ .
\end{array}
\right.
$$
By definition, when $S$ is a structure on $V$, the pair $(V,V^*)$ is a
{\emph{proto-bialgebroid}}. The anchor and bracket of $V$ and of $V^*$ are the
following {\emph{derived brackets}} \cite{yks96} \cite{yks04} 
\cite{yks05} \cite{R} \cite{Vo}:
$$
\begin{array}{llllll}
\vspace{.1cm}
{\mathrm{anchor}} & {\mathrm{of}} & V, & 
a_*(X) \cdot f & = & \{\{X,\mu\},f\} \ , \\
\vspace{.1cm}
{\mathrm{bracket}} & {\mathrm{of}} &  V, &
\mu(X,Y) & = & \{\{X,\mu\},Y\} \ , \\
\vspace{.1cm}
{\mathrm{anchor}} & {\mathrm{of}} &  V^*, &
a^*(\alpha) \cdot f & = & \{\{\alpha,\gamma\},f\} \ , \\
\vspace{.1cm}
{\mathrm{bracket}} & {\mathrm{of}} &  V^*, & 
\gamma(\alpha, \beta) & = & \{\{\alpha,\gamma\},\beta\} \ , \\
\end{array}
$$
for $f \in C^{\infty} (M)$,  \, $X$ and $Y \in \varGamma (V)$, \, $\alpha$
and $\beta \in \varGamma (V^*)$. The quasi-Gerstenhaber brackets on 
$\varGamma(\wedge^\bullet V)$, where $\wedge^\bullet V$ is the
exterior algebra of $V$, and on $\varGamma(\wedge^\bullet
V^*)$, are expressed by the same formulas. They are denoted by
$[~,~]_\mu$ and $[~,~]_\gamma$, respectively. 

The Lie-quasi bialgebroids, quasi-Lie
bialgebroids and Lie bialgebroids are defined as follows:

\noindent $\bullet$ 
$(V,V^*)$ is a 
{\it {Lie-quasi bialgebroid}} if and only if 
$S = \phi + \gamma + \mu $, {\it {i.e.}}, if $\psi = 0$.
Then $V$ is a Lie algebroid,
$\varGamma (\wedge^\bullet V)$ is a Gerstenhaber algebra,
while
$\varGamma ( \wedge^\bullet V^{*})$ is a quasi-Gerstenhaber algebra.

\noindent $\bullet$ 
$(V,V^*)$ is a {\it {quasi-Lie bialgebroid}} 
if and only if $S = \gamma + \mu + \psi$, {\it {i.e.}}, if $\phi =
0$. Then $V^*$ is a Lie algebroid, 
$\varGamma (\wedge^\bullet V^*)$ is a Gerstenhaber algebra,
while
$\varGamma (\wedge^\bullet V)$ is a quasi-Gerstenhaber algebra.

\noindent $\bullet$ 
$(V,V^*)$ is a {\it {Lie bialgebroid}} if and only if 
$S = \gamma + \mu $, {\it {i.e.}}, if $\phi = \psi = 0$.
Then both $V$ and $V^*$ are Lie algebroids, 
and both $\varGamma (\wedge^\bullet V)$
and $\varGamma (\wedge^\bullet V^{*})$ are Gerstenhaber algebras.

The quasi-Gerstenhaber algebras (see \cite{R} \cite{H2} \cite{B1} \cite{B2}) 
are the simplest higher structures beyond the Gerstenhaber algebras themselves;
they correspond to the case where 
all $n$-ary brackets, $\ell_{n}$, vanish for $n\geq 4$.

\medskip

On the Poisson manifold $T^*\varPi V$, we can consider
the Hamiltonian vector field with Hamiltonian
$S \in \mathcal F$, which we denote 
by $d_S = \{ S, . \}$. Because $\{S,S\}=0$, $d_S$ 
is a {\emph{differential}} on the space of
smooth functions on $T^*\varPi V$, {\it{i.e.}}, a derivation of
$\mathcal F$ of degree
$1$ and of square zero.

\medskip

\noindent{\bf Example 1.}
When $V = TM$ and $S = \mu = p_i\xi^i$, then
$\mu(X,Y)$ is the Lie bracket of vector fields $X$ and $Y$, the
corresponding Gerstenhaber bracket on $\varGamma(\wedge^\bullet TM)$ is the
Schouten--Nijenhuis bracket of multivector fields, and the restriction
of $d_S = d_\mu$ to the differential forms on $M$ is the de Rham differential.

\medskip

\noindent{\bf Example 2.} 
When $M$ is a point, then $V= {\mathfrak{g}}$ is a vector space
and a structure $S = \mu +\gamma$ on $V$ is a Lie bialgebra structure on
$({\mathfrak{g}},{\mathfrak{g}}^*)$, also denoted by\break
$S_{\mathfrak{g}} + S_{{\mathfrak{g}}^*}$ in Section \ref{another}, while $d_S = d_\mu +
d_\gamma$ is the Chevalley-Eilenberg cohomology operator of the double
of the Lie bialgebra.
More generally, on $V= {\mathfrak{g}}$,
a structure $S = \mu +\gamma + \phi$, where
$\phi \in \wedge^3 V$, 
is a Lie-quasi bialgebra structure on $({\mathfrak{g}},{\mathfrak{g}}^*)$.
\section{Twisting}\label{2}
We consider a structure $S$ on the vector bundle $V$ that defines a 
proto-bialgebroid structure on $(V,V^*)$, and we shall now study the
{\em{twisting}},
$e^{-\sigma}S$, of $S$ by a function $\sigma$
of shifted bidegree
$(1,-1)$ or $(-1,1)$.
\subsection{Twisting by Poisson or pre-symplectic functions} 
Let $\sigma \in \mathcal F$ be a function
of shifted bidegree
$(1,-1)$ or $(-1,1)$. Since
the right adjoint action, ${\mathrm{ad}}_\sigma = \{ \, .  , \sigma \}$, of an
element $\sigma$ of shifted
degree $0$
is a derivation of degree $0$ of $(\mathcal F, \{~,~\})$, and since, 
for any $a \in \mathcal F$, 
the series $a + \{a,\sigma\} + \frac{1}{2!} \{\{a, \sigma\},\sigma\} 
+ \frac{1}{3!} \{\{\{a, \sigma\},\sigma\}, \sigma\} + \ldots$  
terminates for reasons of bidegrees, the exponential of 
${\mathrm{ad}}_\sigma$ 
is well-defined and is an automorphism of $(\mathcal F,
\{~,~\})$, which, in an abuse of notation, we shall denote by 
$e^{\sigma}$. It follows that,
for any structure $S$, and for any $\sigma$ of shifted degree $0$,
$\{e^{\sigma}S, e^{\sigma}S\} = e^{\sigma}\{S,S\}= 0$, 
and therefore $e^\sigma S$
is also a structure.

\begin{definition}
When $\sigma$
is a function of shifted bidegree $(1,-1)$ or $(-1,1)$, the structure
$e^{-\sigma}S$ is called 
the {\em{twisting}} of $S$ by $\sigma$.
\end{definition}

A function of shifted bidegree $(1,-1)$ is a bivector $\sigma$ on
$V$, 
expressed in local coordinates as
$$\sigma = \frac{1}{2} \sigma^{ab}\theta_a \theta_b \ ,$$
while a function of shifted bidegree $(-1,1)$ is a $2$-form $\tau$ on
$V$, expressed in local coordinates as
$$\tau = \frac{1}{2} \tau_{ab}\xi^a \xi^b \ .$$
We list 
the explicit formulas \cite{R} for the homogeneous components of twisted 
structures.

\noindent $\bullet$ For $\sigma$ of shifted bidegree $(1,-1)$,
let
$e^{- \sigma} S = \phi_{\sigma}
 + \gamma_{\sigma} + \mu_{\sigma} + \psi_{\sigma}$ be the
 decomposition (\ref{structure}) of $e^{-\sigma} S$ as a sum of
 terms of homogeneous bidegrees. Then,

\begin{equation}\label{phi}
\left\lbrace
\begin{array}{ll}
\vspace{.1cm}
&\phi_{\sigma}
= \phi - \{\gamma, \sigma\} + \frac{1}{2} \{\{\mu, \sigma\},\sigma\} -
\frac{1}{6}  
\{\{\{\psi, \sigma\},\sigma\}, \sigma\} \ , \\
\vspace{.1cm}
& \gamma_{\sigma} 
= \gamma - \{\mu, \sigma\} + \frac{1}{2} \{\{\psi, \sigma\},\sigma\} \\
\vspace{.1cm}
&\mu_{\sigma} = \mu - \{\psi,\sigma\} \ , \\
\vspace{.1cm}
& \psi_{\sigma}= \psi \ .
\end{array}
\right.
\end{equation}

%\smallskip

\noindent $\bullet$ For $\tau$ of shifted bidegree $(-1, 1)$,
let
$e^{- \tau} S = \phi_{\tau}
 + \gamma_{\tau} + \mu_{\tau} + \psi_{\tau}$ be the
decomposition (\ref{structure}) of $e^{-\tau} S$ as a sum of
terms of homogeneous bidegrees. Then,

\begin{equation}\label{psi}
\left\lbrace
\begin{array}{ll}
\vspace{.1cm}
\phi_{\tau}= \phi \ , \\
\vspace{.1cm} 
 \gamma_{\tau} 
= \gamma - \{\phi, \tau\} \ , \\
\vspace{.1cm} 
\mu_{\tau} = \mu - \{\gamma,\tau\},
+ \frac{1}{2} \{\{\phi, \tau\},\tau\} \ , \\
\vspace{.1cm} 
\psi_{\tau}
= \psi - \{\mu, \tau\} + \frac{1}{2} \{\{\gamma, \tau\},\tau\} -
\frac{1}{6}  
\{\{\{\phi, \tau\},\tau\}, \tau\} \ . 
 \end{array}
\right.
\end{equation}

\begin{definition} Let $S$ be a structure on $V$.

\noindent (i) 
A function $\sigma$ 
of shifted bidegree
$(1,-1)$ such that $\phi_{\sigma} =0$ 
is called a {\em {Poisson function}} with respect to $S$.

\noindent (ii)
A function $\tau$ 
of shifted bidegree
$(- 1, 1)$ such that 
$\psi_{\tau} =0$
is called a {\em {pre-symplectic function}} with respect to $S$.
\end{definition} 

In view of these definitions, we immediately obtain

\begin{proposition}
Let $S$ be a structure on $V$ and let $\sigma$ (resp., $\tau$) be a
function of shifted bidegree $(1,-1)$ (resp., $(-1,1)$). 

\noindent (i) If $\sigma$ is a Poisson function, the twisted structure
$e^{-\sigma}S$ is a
quasi-Lie\break bialgebroid structure.

\noindent (ii) If $\tau$ is a pre-symplectic function,
the twisted structure $e^{-\tau}S$ is a Lie-quasi bialgebroid
structure. 
\end{proposition}

\subsection{Twisting by Poisson functions}\label{maurer}
It follows from the formula for $\phi_\sigma$ in (\ref{phi}) that a section 
$\sigma$ of $\wedge^2 V$ is a Poisson function with respect to a structure
$S =\phi + \gamma + \mu + \psi$ if and only if
\begin{equation}\label{MC}
\phi - \{\gamma, \sigma\} + \frac{1}{2} \{\{\mu, \sigma\},\sigma\} -
\frac{1}{6}  
\{\{\{\psi, \sigma\},\sigma\}, \sigma\} = 0 \ .
\end{equation}
Equation (\ref{MC}) is called a {\it generalized twisted Maurer-Cartan
equation}, or simply a {\it Maurer-Cartan equation}. 

For any bivector $\sigma$, we set $\sigma^\sharp =i_\alpha \sigma$,
for $\alpha \in \varGamma (V^*)$, where $i$ denotes the interior product.
Whenever $\sigma$ is a Poisson function with respect to $S =\phi +
\gamma + \mu + \psi$, the term of shifted bidegree $(1,0)$ in
$e^{- \sigma}S$,
 $$\gamma_\sigma
= \gamma - \{\mu, \sigma\} + \frac{1}{2} \{\{
\psi, \sigma\},\sigma \} \ ,
$$ 
defines an anchor $a^* + a_* \circ \sigma^\sharp$ and a Lie bracket on 
$\varGamma (V^*)$, as well as a\break Gerstenhaber bracket
on $\varGamma(\wedge^\bullet V^*)$, which we denote by
$[~,~]_{\gamma_\sigma}$, and
a differential $d_{\gamma_\sigma} = \{\gamma_\sigma, . \}$ on    
$\varGamma(\wedge^\bullet V)$.
There is also a bracket, $[~,~]_{\mu_\sigma}$, on
$\varGamma(\wedge^\bullet V)$ defined by the term of shifted bidegree
$(0,1)$,
$\mu_\sigma = \mu - \{\psi, \sigma\}$,
and a derivation of
degree $1$, $d_{\mu_\sigma} = \{\mu_\sigma, . \}$, of 
$\varGamma(\wedge^\bullet V^*)$. Then
$\frac{1}{2} \{\mu_\sigma, \mu_\sigma \} + \{\gamma_\sigma, \psi\} = 0$,
so that $\psi$ measures the defect in the Jacobi
identity for $[~,~]_{\mu_\sigma}$, and $(d_{\mu_\sigma})^2 = [\psi,
\cdot]_{\gamma_\sigma}$.

It appears that the
twisting of Lie bialgebras in the sense of Drinfeld \cite{D2}, as well as its
generalizations to proto-bialgebras \cite{yks92} \cite{BK}
and to proto-bialgebroids
\cite{R}, and the twisting of Poisson structures in the sense 
of {\v S}evera and Weinstein \cite{SW}, and its generalizations 
to structures on Lie algebroids \cite{R} \cite{KL},  
all fit into this general framework, although the meaning of the word
``twisting'' is not quite the same in both instances. In the first
instance, one twists a given structure, in the sense of Definition
\ref{structuredef}, 
on a Lie algebra $\mathfrak g$ 
by an element $\sigma \in \wedge^2 \mathfrak g$ (often denoted by $t$
or $f$), called the ``twist'' \cite{D2} \cite{AK}. For any twist, a
Lie-quasi bialgebra is twisted into a Lie-quasi bialgebra. 
In the second case, it would be more appropriate to speak of ``Poisson
structures with background'': the given structure on the vector bundle
$V$ is of the form $\mu + \psi$, where $\psi$ is a $d_\mu$-closed $3$-form,
and equation (\ref{MC}) which reduces to 
the twisted Poisson condition (\ref{twP}) below is the condition for 
$\sigma \in \varGamma
(\wedge^2 V)$ to twist $\mu + \psi$ into a quasi-Lie bialgebroid
structure.

\smallskip

\noindent {\bf (i) Twist in the sense of Drinfeld}. 
In the case of a twist of a Lie-quasi bialgebra,
one twists a structure $S = \phi
+ \gamma + \mu + 0$ on a Lie algebra $\mathfrak g$
by an arbitrary $\sigma \in \wedge^2 \mathfrak g$ 
into $$
e^{-\sigma}S =
\left( \phi -
\{\gamma, \sigma\} + \frac{1}{2} \{\{\mu, \sigma\},\sigma\} \right) 
+  
\left( \gamma - \{\mu, \sigma\} \right) + \mu + 0 \ ,
$$ 
and one
obtains a ``twisted Lie-quasi bialgebra''. 
The resulting object is a
Lie bialgebra, with $\mu_{\sigma} =\mu$ and $\gamma_{\sigma} = \gamma - 
\{\mu ,\sigma\}$,
if and only if $\sigma$ is a Poisson
function, {\it{i.e.}}, satisfies the condition
$$
\frac{1}{2} [\sigma
,\sigma]_{\mu} + d_{\gamma} \sigma - \phi = 0 \ .
$$
If one
twists a Lie bialgebra ($\psi = \phi = 0$), this condition reduces to
the usual {\it Maurer-Cartan equation},
\begin{equation}\label{MC0}
\frac{1}{2} [\sigma, \sigma]_{\mu} + d_\gamma \sigma = 0 \ .
\end{equation}
If one twists a trivial Lie bialgebra ($\psi = \phi = \gamma = 0$),
the Maurer-Cartan
equation reduces to 
$[\sigma, \sigma]_{\mu} = 0$, {\it{i.e.}}, to the classical Yang-Baxter
equation. 
In fact, for $\sigma = r \in \wedge^2 \mathfrak g$,
$$
- \frac{1}{2} [r,r]_{\mathfrak g} = [r_{12},r_{13}] + [r_{12},r_{23}] +
[r_{13} + r_{23}] \ ,
$$
and the {\it classical Yang-Baxter equation} (CYBE) on a Lie algebra
$\mathfrak g$ is the condition 
$[r_{12},r_{13}] + [r_{12},r_{23}] + [r_{13} + r_{23}] = 0$, for $r
\in \wedge^2 \mathfrak g$.

When $S = \mu$, the necessary and sufficiant
condition for $\mu + \gamma_{\sigma}$ to be a Lie
bialgebra structure on $(\mathfrak g, {\mathfrak g}^*)$ is 
$\{ \mu, \{ \{ \mu, \sigma \}, \sigma \} \} = 0$, the
{\it {generalized classical Yang-Baxter equation}}, which states tht
$[\sigma,\sigma]_\mu$ is ${\rm {ad}}^\mu$-invariant.

In the same way, a Lie-quasi bialgebroid can be twisted by a
bivector, and a Lie bialgebroid is twisted into a Lie bialgebroid if
and only if the bivector satisfies the Maurer-Cartan equation
(\ref{MC0}) (see \cite{LWX} \cite{R} \cite{yks05}).

\smallskip

\noindent {\bf (ii) Twisted Poisson structures}.
If $S$ is a structure on a vector bundle $V$ 
such that $\gamma = 0$ and $\phi = 0$, then 
$\{\mu,\mu\}=0$, 
{\it {i.e.}}, $V$ is a Lie algebroid, and $\psi$ is a $d_{\mu}$-closed section of
$\wedge^3V^*$. In this case, one twists $S= 0 + 0 + \mu + \psi$ into 
$$
e^{-\sigma}S =\left(\frac{1}{2} \{\{\mu, \sigma\},\sigma\} -
\frac{1}{6}  
\{\{\{\psi, \sigma\},\sigma\}, \sigma\} \right) \hspace{3cm}
$$
\vspace{-.4cm}
$$
+ \left(- \{\mu, \sigma\} + \frac{1}{2} \{\{\psi, \sigma\},\sigma\}
\right) + \left(\mu - \{\psi, \sigma\} \right) + \psi \ .
$$ 
Thus, $\sigma$ is a Poisson function 
if and only if 
$$
\{\{\mu, \sigma\},\sigma\} -
\frac{1}{3}  
\{\{\{\psi, \sigma\},\sigma\}, \sigma\} = 0 \ ,
$$
which is the condition
\begin{equation}\label{twP}
\frac{1}{2} [\sigma,\sigma]_\mu= (\wedge^3\sigma^\sharp)\psi \ ,
\end{equation}
{\it {i.e.}}, 
$(\sigma,\psi)$ is a {\it twisted Poisson structure} on the Lie
algebroid $V$. 
When $\sigma$ satisfies the
twisted Poisson condition (\ref{twP}), the resulting object is a 
quasi-Lie bialgebroid. In particular, 
$- \{\mu, \sigma\} + \frac{1}{2} \{\{\psi, \sigma\},\sigma\}$ is a Lie
algebroid bracket on $V^*$.

If, in addition, $\psi = 0$, then $\sigma$ is a Poisson function if and only if 
$$
\{\{\mu, \sigma\},\sigma\} = 0 \ ,
$$
which is the condition
$$[\sigma,\sigma]_\mu= 0 \ ,
$$
{\it {i.e.}}, 
$\sigma$ is a Poisson structure in the usual sense, a section of  
$\wedge^2V$ with Schouten--Nijenhuis square zero. The Poisson case
is also called the {\it triangular} case by extension of the
terminology used in the theory of Lie bialgebras.

\smallskip

\noindent{\bf The twisted differential}.
In the Poisson case ($\gamma = 0$ and $\psi = 0$), the anchor of $V^*$
is $a_* \circ \sigma^\sharp$, 
and the bracket on 
$\varGamma(\wedge^\bullet V^*)$ is  
$\gamma_\sigma = \{\sigma , \mu \}$, 
the Koszul bracket\footnote{The Koszul bracket \cite{Ko}
restricts 
to the bracket of sections of $\varGamma (V^*)$
generalizing the well-known bracket of $1$-forms on a Poisson
manifold. 
The bracket of $1$-forms on 
symplectic manifolds was introduced in the book of Abraham and
Marsden (1967). For Poisson manifolds, 
it was discovered independently in the 1980's by 
several authors -- Gelfand and Dorfman,
Fuchssteiner, Magri and Morosi, Daletskii --, and Weinstein \cite{CDW}
has shown that it is a Lie algebroid bracket.}. 
The corresponding differential
on $\varGamma(\wedge^\bullet
V)$ is the 
{\it Lichnerowicz-Poisson differential} \cite{Li}, 
$d_\sigma =
\{\{\sigma,\mu\}, . \} = [\sigma, . ]_{\mu}$, while the differential on 
$\varGamma(\wedge^\bullet V^*)$ is the Lie algebroid cohomology operator
$d_\mu= \{ \mu, . \}$. The pair $(V,V^*)$ is a Lie algebroid.

In the twisted Poisson case, $\gamma_\sigma = - \{\mu,
\sigma\} + \frac{1}{2} \{\{\psi, \sigma\},\sigma \}$ restricts to the 
Lie algebroid 
bracket on sections of $V^*$ defined by \v Severa and Weinstein 
\cite{SW},
and the corresponding 
differential on $\varGamma(\wedge^\bullet V)$
is the {\it twisted Poisson differential},
$d_{\sigma} + i_{\psi^{(2)}}$,
where 
$\psi^{(2)} = \frac{1}{2}\{\{\psi, \sigma\}, \sigma\} = 
(\wedge^2 \sigma^{\sharp}) \psi$, while 
the derivation $\{\mu_\sigma, . \}$ is the derivation 
$d_\mu + i_{\psi^{(1)}}$, where $\psi^{(1)} = \{\psi, \sigma\} =
\sigma^{\sharp} \psi$ (see \cite{SW} \cite{R} \cite{KL}).
The pair $(V,V^*)$ is then a quasi-Lie bialgebroid.

\subsection{Twisting by pre-symplectic functions}
It follows from formula (\ref{psi}) that a
section $\tau$ of $\wedge^2 V^*$ is a pre-symplectic function with respect to a structure
$S =\phi + \gamma + \mu + \psi$ if and only if
\begin{equation}\label{MC2}
\psi - \{\mu, \tau\} + \frac{1}{2} \{\{\gamma, \tau\},\tau\} -
\frac{1}{6}  
\{\{\{\phi, \tau\},\tau\}, \tau\} = 0 \ .
\end{equation}
Equation (\ref{MC2}) is dual to (\ref{MC}) and it 
is also called a {\it generalized twisted Maurer-Cartan
equation} or again simply a {\it{Maurer-Cartan equation}}. Pre-symplectic functions generalize pre-symplectic
structures on manifolds
as well as their twisted versions.

If $\gamma = \phi = 0$, then 
$\{\mu,\mu\}=0$,
{\it {i.e.}}, $V$ is a Lie algebroid, 
and $\psi$ is a $d_{\mu}$-closed section of
$\wedge^3V^*$. In this case, $\tau$ is pre-symplectic if and only if 
the pair $(\tau, \psi)$ satisfies the twisted pre-symplectic condition,
$$
\psi - \{\mu, \tau\} = 0 \ ,
$$ 
which is the condition,
$d_\mu \tau = \psi$, {\it {i.e.}}, $(\tau, \psi)$ is a {\it {twisted pre-symplectic structure}} on
the Lie algebroid $V$.
(See \cite{SW} and see \cite{S} for an example of a twisted
symplectic structure
arising in the theory of the lattices of Neumann oscillators.)

If, in particular, $\gamma = \phi = \psi = 0$, 
then $\{\mu,\mu\}=0$
and $V$ is a Lie algebroid. 
In this case, $\tau$ is pre-symplectic if and only if $\tau$
satisfies
the pre-symplectic condition,
$$\{\mu, \tau\} = 0 \ ,
$$
which is the condition,
$d_\mu \tau = 0$, {\it {i.e.}}, $\tau$ is a $d_{\mu}$-closed section of
$\wedge^2V^*$, the {\it {pre-symplectic}} case.
\section{The graphs of Poisson and of pre-symplectic
  functions}\label{3}

\subsection{Courant algebroids, the Courant algebroid $V \oplus V^*$}
A {\it Loday algebra} (called {\it Leibniz algebra} by Loday \cite{Lo}) is equipped
with a bracket (in general non skew-symmetric) satisfying the 
Jacobi identity in the form $[u,[v,w]] = [[u,v],w] + [v,[u,w]]$.
We give the definition of Courant algebroids in \cite{yks05} which is
equivalent to the original definition of Courant and Weinstein 
\cite{CW} \cite{C}.

A {\it Courant algebroid}
is a vector bundle  
$E \to M$, equipped with a vector bundle morphism, $a_E : E \to TM$, called the
{\it anchor}, a fiber-wise non-degenerate symmetric bilinear form $(~,~)$, and
a bracket, $[~,~] : \varGamma (E) \times \varGamma (E) \to \varGamma (E)$,
called the {\it Dorfman-Courant bracket},
such that  

\noindent $\bullet$ 
$\varGamma (E)$ is a {\it Loday algebra},

\noindent $\bullet$ for all $x$, $u$, $v \in \varGamma (E)$,
$$
a_{E}(x)\cdot (u,v) = (x\, , \, [u,v] + [v,u])
= ([x,u] \, , \, v) + (u \, , \, [x,v]) \ .
$$

A sub-bundle, $F \subset E$, is called a {\it Dirac sub-bundle} if

\noindent $\bullet$ $F$ is maximally isotropic,

\noindent $\bullet$ $\varGamma (F)$ is closed under the bracket.

\medskip

When $S$ is a structure on $V$, the vector bundle 
$E = V \oplus V^{*}$ with the canonical scalar product,
$$
(u,v) = \{u,v\} \ ,
$$
and bracket 
\begin{equation} \label{courant} 
[u,v]_S = \{\{u, S\},v\} \ ,
\end{equation}
for $u$, $v \in \varGamma (V \oplus V^*)$,
is a Courant algebroid \cite{R} \cite{Vo} \cite{yks05}, called the
{\it double} of $V$.
\begin{lemma} Let $S$ be a structure on $V$.

\noindent (i) The function $\sigma \in \varGamma(\wedge^2 V)$ 
is a Poisson function with respect to $S$ if and only if $V^*$ is
a Dirac sub-bundle of $(V \oplus V^*, [~,~]_{e^{-\sigma}S})$.

\noindent (ii) The function $\tau \in \varGamma(\wedge^2 V^*)$ is a 
pre-symplectic function with respect to $S$ if and only if $V$ is
a Dirac sub-bundle of $(V \oplus V^*, [~,~]_{e^{-\tau}S})$.
\end{lemma}\label{dirac}
\begin{proof} Part (i) (resp., (ii))
follows from the computation of the bidegrees of the homogeneous
terms in $[u,v]_{e^{-\sigma}S}$ (resp., $[u,v]_{e^{-\tau}S}$) 
for $u, v \in \varGamma (V)$ (resp., $u, v \in \varGamma (V^*)$).
\hfill{$\square$}
\end{proof}

\subsection{Graphs as Dirac structures} 

Theorem \ref{graph} below generalizes the characterization of the graphs
of Poisson, quasi-Poisson and pre-symplectic structures
in \cite{LWX} and \cite{SW},
and that of twisted pre-symplectic structures in
\cite{AX} and \cite{BC}. The statement of this theorem 
can be found in Remark 4.2 in \cite{R} (cf. also Prop. 5 in \cite{BCS}), 
and the proof given here is also
due to Roytenberg \cite{Roytenberg}. Both theorems in this section
have been proved by Terashima \cite{T}.

\begin{theorem} \label{graph}
Let $S$ be a structure on $V$.

\noindent(i) A section $\sigma$ of $\wedge^2 V$ 
is a Poisson function with respect to $S$
if and only if its graph in the Courant algebroid 
$(V \oplus V^*, [~,~]_S)$ is a Dirac sub-bundle.

\noindent(ii) A section $\tau$ of $\wedge^2 V^*$ is a pre-symplectic function 
with respect to $S$ if and only  if its graph in the 
Courant algebroid $(V \oplus V^*, [~,~]_S)$ is a Dirac sub-bundle.
\end{theorem}
\begin{proof}
We need only prove (ii), since the proof of (i) is entirely similar.
We shall denote by $\tau^\flat$ the vector bundle morphism from $V$ to $V^*$
induced by $\tau \in \varGamma(\wedge^2 V^*)$, such that 
$\tau^{\flat}X = - i_X \tau$, for $X \in V$, as well as the associated
map on sections of $V$.
By the graph of $\tau$, we mean the graph of $\tau^\flat$.
Since $\tau^\flat (X) = \{X , \tau \}$, for all $X \in \varGamma (V)$. 
and since, for reasons of bidegree,
$e^{\tau} X = X + \{ X , \tau\}$, it follows that 
\begin{equation} \label{exptau} 
{\rm {Graph}}(\tau)= e^{\tau} V \ .
\end{equation}

Since $e^\tau$ is an automorphism of $(\mathcal F, \{~,~\})$, it is an
isomorphism from $(V \oplus V^*, [~,~]_{e^{-\tau}S})$ to 
$(V \oplus V^*, [~,~]_S)$. Thus $e^\tau V$ is a Dirac sub-bundle of 
$(V \oplus V^*, [~,~]_{S})$ if and only if $V$ is a Dirac
sub-bundle of $(V \oplus V^*, [~,~]_{e^{-\tau}S})$.
Thus (ii) follows from (\ref{exptau}) and Lemma \ref{dirac} (ii).
\hfill{$\square$}
\end{proof}
\begin{theorem} Let $S = \phi + \gamma + \mu + \psi$ be a structure on $V$.

\noindent(i) Let $\sigma$ be a Poisson function with respect to $S$. The
projection ${\mathrm{Graph}}(\sigma) \to \varGamma (V^*)$ 
is a morphism of Lie algebroids
when ${\mathrm{Graph}}(\sigma)$ is equipped with the Lie bracket
induced from the Dorfman-Courant bracket $[~,~]_S$
and  
$\varGamma (V^*)$ is equipped with the Lie bracket 
$
\gamma_{\sigma} = \gamma - \{\mu, \sigma\} + \frac{1}{2} \{\{\psi,
\sigma\},\sigma\}$.                                        

\noindent(ii) Let $\tau$ be a pre-symplectic 
function with respect to $S$. The
projection ${\mathrm{Graph}}(\tau) \to \varGamma(V)$ 
is a morphism of Lie algebroids
when ${\mathrm{Graph}}(\tau)$ is equipped with the Lie bracket
induced from the Dorfman-Courant bracket $[~,~]_S$
and  
$\varGamma(V)$ is equipped with the Lie bracket 
$
\mu_{\tau} = \mu - \{\gamma, \tau\} + \frac{1}{2} \{\{\phi,
\tau\},\tau\}$.                                        
\end{theorem}
\begin{proof}
We need only prove (ii), since the proof of (i) is entirely similar.
For any $\tau \in \varGamma (\wedge^2 V^*)$, $X$ and $Y
\in \varGamma (V)$, \, 
$[e^\tau X, e^\tau Y]_S
 = e^\tau [X,Y]_{e^{-\tau}S}$.
If $\tau$ is a
pre-symplectic function with respect to $S$, then 
$[e^\tau X, e^\tau Y]_S
 = e^\tau [X,Y]_{\mu_\tau} = [X,Y]_{\mu_\tau} + \{ [X,Y]_{\mu_\tau},
 \tau\}$, 
whose $V$-component is $[X,Y]_{\mu_\tau}$. 
\hfill{$\square$}
\end{proof}

\section{Symplectic functions}\label{4}
Let us now assume that 
$\sigma \in \varGamma( \wedge^2 V)$ 
is non-degenerate, {\it{i.e.}},
the map $\sigma^{\sharp} : V^* \to~V$ defined by $\sigma^{\sharp} \alpha =
i_{\alpha }\sigma$, for $\alpha \in \varGamma (V^*)$, 
is invertible. Set $\tau^{\flat} =
(\sigma^{\sharp})^{-1}$,
and let $\tau \in \varGamma(\wedge^2 V^*)$ be such that 
$\tau^{\flat}X = - i_X \tau$, for $X \in \varGamma (V)$. 
We say that $\tau \in \varGamma( \wedge^2 V^*)$ and $\sigma \in \varGamma(
\wedge^2 V)$ are 
 {\it inverses} of one another.
A non-degenerate pre-symplectic function is called {\it symplectic}.

\subsection{``Non-degenerate Poisson'' is equivalent to ``symplectic''}
Many classical results are corollaries of the general
theorem which we state and prove in this section. 
Recall that 
$\xi^a \theta_b = - \theta_b \xi^a$,
$\{\xi^a ,\theta_b \} =  \delta^a_b = \{\theta_b , \xi^a \}$
and, for $u, v, w \in \mathcal F$,
$$
\{u,vw\}= \{u,v\} w + (-1)^{|u||v|} v \{u,w\} \ , 
$$
$$
\{uv, w\}= u \{v, w\} + (-1)^{|v||w|} \{u,w\} v \ ,
$$
where $|u|$ is the degree of $u$, and
$$ \{u,\{v,w\}\}=
\{\{u,v\}, w \}+ (-1)^{\|u\| \|v\|} \{v, \{u,w\}\} \ ,
$$
$$
\{\{u,v\}, w\}= \{u,
\{v, w\}\} + (-1)^{\|v\| \|w\|} \{\{u,w\}, v\}  \ ,
$$
where $\|u\|$ is the shifted degree of $u$.
The proof of the theorem depends on the following lemma. 
\begin{lemma}\label{lemma}
Assume that $\sigma \in \varGamma(\wedge^2V)$ 
is non-degenerate and that its inverse is $\tau$. Then

\noindent (i) $\{\sigma, \tau \} = - \{\tau, \sigma \}  
= {\mathrm{Id}}_V$.

\noindent (ii) If $S$ is of shifted bidegree $(p,q)$,
then
\begin{equation}\label{p-q}
\{\{\sigma,\tau\}, S\} = (q-p) S \ .
\end{equation}
\end{lemma}
\begin{proof}
This lemma is proved by straightforward computations, 
using the equality ${\mathrm{Id}}_V = \xi^a \theta_a $. 
\hfill$\square$
\end{proof}

\begin{theorem}\label{symplectic} Let $S$ be a structure on $V$.
Let $\sigma \in \varGamma( \wedge^2 V)$ be a non-degenerate bivector
with inverse $\tau \in \varGamma (\wedge^2 V^*)$. 
Then $\sigma$ is a Poisson function with respect to $S$ if and
only if $ - \tau$ is a symplectic function with respect to $S$.
\end{theorem}
\begin{proof}
Lemma \ref{lemma}(ii) applied in the cases $(p,q) = (2,-1)$, $(1,0)$, $(0,1)$ and
$(-1,2)$, and repeated applications of the Jacobi identity yield 
the following computations.  
Let $\mu$ be of shifted bidegree $(0,1)$. From 
$$\{\{\mu,
\tau\},\sigma\} = \{\mu,
\{\tau,\sigma\} \} + \{\{\mu, \sigma\},\tau\} =  \mu +
 \{\{\mu, \sigma\},\tau\} \ ,
$$ 
we obtain 
$$\{\{\{\mu,
\tau\},\sigma\},\sigma\} =
\{\mu,\sigma\} + \{\{\{\mu, \sigma \}, \tau\},\sigma\}
$$
$$
= 
\{\mu,\sigma\} + \{\{\mu, \sigma \}, \{\tau, \sigma\}\} + \{\{\{\mu, \sigma\},\sigma\},\tau\}
=
\{\{\{\mu,
\sigma\},\sigma\},\tau\} \ .
$$
Whence
$$
\{\{\{\{\mu, \tau\},\sigma\}, \sigma\}, \sigma\}=
\{\{\{\{\mu, \sigma\},\sigma\},\tau\},\sigma\} = 
\{\{\{\mu, \sigma\},\sigma\}, \{\tau,\sigma\}\}  
$$
$$
= - \, 3 \, \{\{\mu, \sigma\},\sigma\} \ .
$$
Similarly, if $\gamma$ is of shifted bidegree $(1,0)$, 
$$\{\{\{\{\{\gamma, \tau\}, \tau\},\sigma\}, \sigma\}, \sigma\}= 12 \,
\{\gamma,\sigma\} \ .
$$
If $\phi$ is of shifted bidegree $(2,-1)$,
$$
\{\{\{\{\{\{\phi, \tau\}, \tau\}, \tau\},\sigma \}, \sigma\}, \sigma\}
= - \, 36 \, \phi \ .
$$

Let  $S= \phi + \gamma + \mu + \psi$.
The term of shifted bidegree $(-1,2)$ in
$e^{-\tau} S$ is 
$$\psi_{\tau}= 
\psi - \{\mu, \tau\} + \frac{1}{2} \{\{\gamma, \tau\},\tau\} -
\frac{1}{6}  
\{\{\{\phi, \tau\},\tau\}, \tau\}\ ,
$$ 
and the term of shifted bidegree $(2,- 1)$ in
$e^{-\sigma} S$ is 
$$\phi_{\sigma} = 
\phi - \{\gamma, \sigma\} + \frac{1}{2} \{\{\mu, \sigma\},\sigma\} -
\frac{1}{6}  
\{\{\{\psi, \sigma\},\sigma\}, \sigma\} \ .
$$
The preceding equalities 
and analogous results for other 
iterated brackets, reversing the roles of
$\sigma$ and $\tau$, yield the equalities:
$$
\{\{\{ \psi_{\tau}, \sigma\},\sigma\},\sigma\}
=  6 \, \phi_{- \sigma}
$$
and
$$
\{\{\{ \phi_{\sigma}, \tau\},\tau\},\tau\}
= 6 \, \psi_{- \tau} \ .
$$
Therefore $\psi_{\tau} = 0$ implies $\phi_{-\sigma} = 0$, and conversely. 
\hfill{$\square$}
\end{proof}

The method of proof 
used above 
in the general case can be applied to give one-line proofs of
some well-known results.

\noindent
$\bullet$
For the case of non-degenerate Poisson structures,
the proof reduces to $\{\mu, \tau\}=0$ implies that
$\{\{\{\{\mu,
\tau\},\sigma\}, \sigma\}, \sigma\}
= 0$, which implies that\break 
$\{ \{ \mu, \sigma \},\sigma \} = 0$, and a similar argument applies to the converse. 
This simple argument proves the classical result: non-degenerate
closed $2$-forms are in one-to-one correspondence with non-degenerate
Poisson bivectors.

\noindent$\bullet$ 
For the case of non-degenerate twisted Poisson structures (see 
Section \ref{maurer} (ii)), the proof reduces
to $\{\mu, \tau\} = - \psi$ implies that $\{\{\{\{\mu,
\tau\},\sigma\}, \sigma\}, \sigma \} = - \{ \{\{\psi,
\sigma\},\sigma\}, \sigma\}$, which implies  
that $\{\{\mu,
\sigma\},\sigma\}
= \frac{1}{3}
\{\{\{\psi,
\sigma\},\sigma\}, \sigma\}$, and a similar argument for the converse.
Thus $d_\mu \tau = - \, \psi$ implies $\frac{1}{2}[\sigma,\sigma]_\mu =
(\wedge^3\sigma^\sharp)\psi$ and conversely.
This constitutes a direct proof of the 
following corollary 
of Theorem \ref{symplectic} (see \cite{SW} \cite{AX} \cite{KY}).

\begin{corollary}
(i) A non-degenerate bivector on a Lie algebroid
defines a twisted Poisson structure 
if and only if its inverse is a twisted symplectic\break $2$-form.

\noindent (ii) The leaves of a twisted Poisson manifold are twisted
symplectic manifolds.
\end{corollary}

It follows from this corollary 
that, in the case of Lie algebras, considered to be Lie
algebroids over a point,
a non-degenerate $r \in \wedge^2 {\mathfrak{g}}$ is a solution 
of the {\it{twisted classical Yang-Baxter equation}},
generalizing the classical
Yang-Baxter equation (see Section \ref{maurer}),
$$ \frac{1}{2} [r,r]_{\mathfrak g} = (\wedge^3 r^{\sharp}) \psi \ ,$$
where $\psi$ is a $d_{{\mathfrak{g}}}$-closed 
$3$-form on the Lie algebra ${\mathfrak{g}}$, if and only if its inverse
is a non-degenerate $2$-form $\tau$ satisfying the twisted closure
condition,
$d_{{\mathfrak{g}}} \tau =
- \psi$.
Here $d_{\mathfrak g}$ is the Chevalley-Eilenberg cohomology operator of
$\mathfrak g$ and
the bracket, $[~,~]_{\mathfrak g}$, is the
algebraic Schouten bracket on $\wedge^\bullet \mathfrak g$.

Recall that a Lie algebra is called {\it
quasi-Frobenius} if it possesses a non-degenerate $2$-cocycle. 
Thus, we recover in particular the well-known correspondence \cite{Sto}
\cite{GG} \cite{HY}
between 
non-degenerate triangular $r$-matrices, {\it{i.e.}}, skew-symmetric
solutions of the classical Yang-Baxter equation,  and quasi-Frobenius
structures.
\begin{corollary}
A non-degenerate bivector in $\wedge^2 {\mathfrak{g}}$ is a 
solution of the classical Yang-Baxter equation if and only if its
inverse defines a quasi-Frobenius structure on ${\mathfrak{g}}$.
\end{corollary}

\subsection{Regular twisted Poisson structures}
We summarize a result from \cite{KY} which can now be considered to be
a corollary of Theorem \ref{symplectic}.
Let $A$ be a vector bundle with a bivector $\pi \in \varGamma (\wedge^2
A)$ such that $\pi^{\sharp}$ is of constant rank.
Let $B$ be the image of $\pi^{\sharp}$. Then $B$ is a Lie
sub-algebroid of $A$ and, because $\pi $
is skew-symmetric, $\pi^{\sharp}$ defines an isomorphism,
$\pi^{\sharp}_B 
: B^* \to B$, where $B^*= A^*/\ker
{\pi^{\sharp}}$ is the dual of $B$.
Then the inverse of $\pi^{\sharp}_B $ defines a non-degenerate 
$2$-form on $B$,
$\omega_B \in \varGamma (\wedge^2 B^*)$, by $(\pi^{\sharp}_B)^{-1} X = - i_X
\omega_B$, for $X \in \varGamma (B)$.

Assume that the vector bundle, $A$, is in fact a Lie algebroid. Let
$\psi$ be a $d_A$-closed
$3$-form on $A$, and let $\psi_B$ denote 
the pull-back of $\psi$ under the canonical injection
$\iota_B : B \hookrightarrow A$. Then

\begin{proposition}
Under the preceding assumptions,
$(A, \pi, \psi)$ is a Lie algebroid with a regular twisted Poisson
structure if and only if $(B, \omega_B, \psi_B)$ is a Lie algebroid
with a twisted
symplectic structure, {\it {i.e.}}, if and only if
$d_B \omega = - \psi_B$.
\end{proposition}

This proposition constitutes a linearization of the twisted Poisson
condition, and can be applied in particular to the case of Lie
algebras \cite{KY}.

\section{Another type of Poisson function:
Lie algebra actions on manifolds}\label{another}

In this section, we consider the twisting of various structures
involving the action of a Lie algebra on a manifold.

\subsection{Structures on $TM \times {\mathfrak g}^*$}
Let $\mathfrak{g}$ be a Lie algebra, and let $M$ be a manifold.
We consider the vector\break 
bundle $V = TM \times \mathfrak{g}^{*}$ over $M$ 
which is, by definition, $TM \mathop{\oplus}\limits_{M} (M \times
\mathfrak{g}^{*}) \to~M$. 
We introduce local coordinates on $T^{*}\varPi V$,
$(x^{i}$, $\xi^{i}$, $e_A$, $p_{i}$, $\theta_{i}$, $\epsilon^{A})$, 
where\break $i= 1, \ldots, \dim M$, and $A= 1,
\ldots, \dim {\mathfrak {g}}$,  
with the following bidegrees,
$$
\begin{array}{ccccccc}
x^{i} & \xi^{i} & e_A & p_{i} & \theta_{i} & \epsilon^{A} & \\
(0,0) & (0,1) & (0,1) & (1,1) & (1,0) & (1,0) & \mathrm{bidegree}\\
(-1,-1) & (-1,0) & (-1,0) & (0,0) & (0,-1) & (0,-1) & 
\mathrm{shifted \, \, bidegree}\\
\end{array}
$$ 
satisfying
$$
\{x^i,p_j\}=\delta^i_j \, , \quad  \,
\{\xi^i,\theta_j\}=\delta^i_j \, , \quad \,
\{e_A,\epsilon^B\}=\delta^B_A \, .
$$

Let 
$$
S_{\mathfrak g} =\frac{1}{2} C^{D}_{AB} \epsilon^{A} \epsilon^{B} e_{D}
$$ 
be the function on $T^{*}\varPi V$ of shifted bidegree $(1,0)$
defining the Lie bracket of
$\mathfrak{g}$, and let 
$$
S_{M} = p_{i} \xi^{i}
$$
be the function on $T^{*}\varPi V$ 
of shifted bidegree $(0,1)$ which defines the Schouten--Nijenhuis bracket of
multivectors on $M$. Then
\begin{equation}\label{liealg}
[u,v]_{\mathfrak{g}} = \{\{ u, S_{\mathfrak{g}} \},v\} \ ,
\end{equation}
for all $u$, $v \in \mathfrak g$, and 
\begin{equation}\label{liebracket}
[X,Y]_{M} = \{\{ X,S_{M}\} , Y\} \ ,
\end{equation}
for all $X$, $Y \in \varGamma (TM)$.
It is easy to show that $S_{\mathfrak {g}} + S_{M}$ is a structure on $V$.

More generally, consider the following functions 
on $T^{*}\varPi V$ of shifted
bidegree $(-1,2)$,
a $3$-form $\Psi_M$ on $M$,
$$
\Psi_{M} = \frac{1}{6} \Psi_{ijk}\xi^{i}\xi^{j}\xi^{k} \ ,
$$
and a $3$-form $\Psi_{\mathfrak{g}}$ on $\mathfrak{g}^{*}$,
$$
\Psi_{\mathfrak{g}} = \frac{1}{6} \Psi^{ABC} e_{A} e_{B} e_{C} \ .
$$
Then 
$S_{\mathfrak{g}} + S_{M} + (\Psi_{\mathfrak{g}} + \Psi_{M})$ is a
structure on $V$ if and only if

\noindent $\bullet \{ S_{M} ,\Psi_{M} \} = 0$, {\it{i.e.}}, 
$\Psi_{M}$ is a closed $3$-form on $M$, and

\noindent$
\bullet \{ S_{\mathfrak{g}} ,\Psi_{\mathfrak{g}} \} = 0$,
{\it{i.e.}},
$\Psi_{\mathfrak{g}}$ is a $0$-cocycle on $\mathfrak{g}$ with
values in $\wedge^{3}\mathfrak{g}$.

More generally still, we can, in addition, introduce a 
function  on $T^{*}\varPi V$
of shifted bidegree $(0,1)$ which defines a  
bracket on $\mathfrak{g}^{*}$,
$$
S_{\mathfrak{g^{*}}} = \frac{1}{2} \Gamma^{AB}_{C} e_{A} e_{B} \epsilon^{C} 
\ .
$$
Then 
$S= S_{\mathfrak{g}} + (S_{\mathfrak{g}^{*}} + S_{M})+
(\Psi_{\mathfrak{g}} + \Psi_{M})$,
a sum of terms of shifted bidegrees $(1,0)$, $(0,1)$ and $(-1,2)$,
respectively, is a structure on $V$ if and only if

\noindent
$\bullet$ $\{ S_{M} ,\Psi_{M} \} = 0$, {\it{i.e.}}, 
$\Psi_{M}$ is a closed $3$-form on $M$, and

\noindent
$\bullet$ 
$\{ S_{\mathfrak g} + S_{{\mathfrak g}^*} +
\Psi_{\mathfrak g} , S_{\mathfrak g} + S_{{\mathfrak g}^*} +
\Psi_{\mathfrak g} \} = 0$, the condition that 
$( \mathfrak{g}, \mathfrak{g}^{*} )$ be a {\it
  Lie-quasi bialgebra}. 

Let us assume that these conditions are satisfied. By 
what function can we twist the structure
$S_{\mathfrak{g}} + (S_{\mathfrak{g}^{*}} + S_{M})+
(\Psi_{\mathfrak{g}} + \Psi_{M})$? We can twist it
by any function of shifted bidegree $(1,-1)$.
Therefore we can choose 
$$
{\rho = \rho^{i}_{A} \epsilon^{A}\theta_{i}} \ ,
$$
and twist $S$ by $\rho$, and/or we can twist $S$ by the bivector 
$$
{\pi = \frac{1}{2} \pi^{ij} \theta_{i} \theta_{j}} \ .
$$
We shall now prove, following Terashima \cite{T}, that twisting by $\rho +\pi$  
provides a natural and unified way of determining 
the Lie algebroid structures
discovered by Lu \cite{L} and by Bursztyn, Crainic and \v Severa
\cite{BC} \cite{BCS}. This method yields an immediate proof of the
fact that these are indeed Lie algebroid structures.

\subsection{Twisting by a Lie algebra action}
Let us first determine the meaning of the condition that $\rho$ be a
Poisson function with respect to $S = S_\mathfrak g + S_M$.
We remark that $\rho$, considered either as a function on $T^*\varPi V$ or
as a map from $\mathfrak g$ to $\varGamma (TM)$, satisfies, for all
$u\in \mathfrak{g}$, 
$$
\{ \rho , u\} = \rho (u) \ .
$$
Computing the terms of shifted bidegrees $(2,-1)$, $(1,0)$ and $(0,1)$ 
of 
the twisted structure, $e^{-\rho}S$, we obtain 
$$
e^{-\rho} (S_{\mathfrak{g}} + S_{M} ) = \left( - \{ S_{\mathfrak{g}} ,\rho \} + \frac{1}{2} \{ \{ S_{M},\rho \}
,\rho \}\right) + (S_{\mathfrak{g}} - \{
S_{M},\rho\}) + S_{M} \ .
$$
Therefore $\rho$ is a Poisson function with respect to $S=S_{\mathfrak g} + S_M$
if and only if
\begin{equation}\label{representation}
- \{ S_{\mathfrak g},\rho\} + \frac{1}{2} \{ \{ S_{M} ,\rho \} ,\rho \}
= 0 \ .
\end{equation}
\begin{lemma}\label{action}
The function $\rho$ is a Poisson function with respect to
$S_{\mathfrak g} + S_M$ if and
  only if it is a Lie algebra action of $\mathfrak{g}$ on $M$.
\end{lemma}
\begin{proof} The proof of the fact that 
relation (\ref{representation}) is equivalent to
$$
\rho([u,v]_{\mathfrak{g}}) =
[\rho(u),\rho (v)]_{M}
$$
for all $u$, $v \in \mathfrak{g}$,
depends on formulas (\ref{liealg}) and
  (\ref{liebracket}),
the Jacobi identity and the
vanishing of all brackets of the form 
$\{e_A, \theta_i \}$ and $\{\epsilon^A, \theta_i \}$, whence
$$
\rho([u,v]_{\mathfrak g}) = \{\{\{S_{\mathfrak g}, \rho \},u \}, v\}
$$
and 
$$ \quad \quad  \quad\quad  \quad \quad \quad \quad  [\rho(u),\rho(v)]_M= \frac{1}{2} \{\{\{\{ S_M,\rho\},\rho\},u\},v\}
\ . \quad \quad \quad \quad \quad
\quad \quad 
{\square} 
$$ 
\end{proof}

\subsection{Introducing additional twisting by a bivector}
Let us now twist $S = S_{\mathfrak{g}} + (S_{\mathfrak{g}^{*}} + S_{M})+
(\Psi_{\mathfrak{g}} + \Psi_{M})$ by 
$$
{\sigma =\pi +\rho} \ .
$$
We first observe that the brackets 
$\{ \pi ,\rho \} $, $\{ S_{\mathfrak{g}},\pi\}$,
$\{S_{\mathfrak{g}*},\pi\} $,
$\{\{ S_{\mathfrak{g}*}, \rho \} ,\pi\} $ and $\{ \Psi_{\mathfrak{g}}
,\pi \} $ vanish. Computing the term of shifted bidegree $(2,-1)$ in 
$e^{-(\pi + \rho)}S$, 
we see that 
$\pi +\rho$ is a Poisson function with respect to $S$ if and only if
$$
\begin{array}{llll}
\vspace{.1cm}
- \{ S_{\mathfrak{g}}, \rho\} + 
\frac{1}{2} \{\{ S_{\mathfrak{g}*} ,\rho\},\rho\}
& + & \frac{1}{2} \{ \{ S_{M} , \pi + \rho \} , \pi + \rho \} \hspace{2.7cm} \\
\vspace{.1cm}
& - & \frac{1}{6} \{\{\{ \Psi_{\mathfrak{g}}+ \Psi_{M} , \pi + \rho \} , \pi
+ \rho \} , \pi + \rho\} = 0 \ .
\end{array}
$$
The computation of the several terms in this generalized
  twisted Maurer-Cartan equation yields 
\begin{proposition}\label{pi+rho} 
The function
$\pi + \rho$ is a Poisson function 
with respect to $S = S_{\mathfrak{g}} + (S_{\mathfrak{g}^{*}} + S_{M})+
(\Psi_{\mathfrak{g}} + \Psi_{M})$ 
if and only if the following four
conditions are satisfied:
$$
\{\{\{ \Psi_{M} ,\rho \} , \rho\} ,\rho \} = 0 \ , \leqno(A)
$$
\vspace{-.7cm}
$$
- \{ S_{\mathfrak{g}},\rho \} + \frac{1}{2} \{\{ S_{M} , \rho\}
,\rho\} - \frac{1}{2} \{\{\{ \Psi_{M} ,\rho \} , \rho \} ,\pi \} 
= 0 \ , \leqno(B)
$$
\vspace{-.6cm}
$$
\{ \{ S_{M} ,\pi \} ,\rho \} + \frac{1}{2} \{ \{
S_{\mathfrak{g}*},\rho \} , \rho \} - \frac{1}{2} \{\{\{ \Psi_{M}
,\rho \} , \pi\} ,\pi \} = 0 \ , \leqno(C)
$$
\vspace{-.6cm}
$$
\{ \{ S_{M},\pi\},\pi\} - \frac{1}{3} \{\{\{
\Psi_{\mathfrak{g}},\rho\} , \rho\} , \rho \}  
- \frac{1}{3} \{\{\{ \Psi_{M} ,\pi \} , \pi\} ,\pi \} = 0 \ . \leqno(D)
$$
\end{proposition}

Condition ($A$) is the relation
$i_{\rho(u) \wedge \rho(v) \wedge \rho(w)} \Psi_M = 0$,
for all $u, v, w \in \mathfrak g$, which means that $\Psi_M$ is in the
kernel of $\wedge^3 \rho^*$, where $\rho^*$ is the dual of $\rho$.

Condition ($B$) is the relation 
\begin{equation}\label{rep}
\rho([u,v]_{\mathfrak{g}}) -
[\rho(u),\rho (v)]_{M} = \pi^\sharp(i_{\rho(u) \wedge \rho(v) } \Psi_M ) \ ,
\end{equation}
for all $u, v \in \mathfrak g$. This is proved by the same
computations as in
Lemma \ref{action}. Thus 
($B$) expresses the fact that 
$\rho$ is a {\it{twisted action}} of
$\mathfrak{g}$ on $M$.

Condition ($C$) is the relation 
\begin{equation}\label{infpoisson}
\mathcal{L}_{\rho (u)} \pi =  - 
(\wedge^{2} \rho) (\gamma (u)) + (\wedge^2 \pi^\sharp)(i_{\rho(u) } \Psi_M ) \ ,
\end{equation}
for all $u \in \mathfrak g$, 
where $\gamma : \mathfrak{g} \to \wedge^{2} \mathfrak{g}$ is 
$S_{\mathfrak{g}*}$ viewed as a cobracket on $\mathfrak{g}$.
In fact, 
$$\{\{\{ S_{M},\pi \} ,\rho \} ,u \} = \{\{\{ \rho ,u\}
,S_{M}\} ,\pi\}
=[ \{ \rho ,u \} ,\pi ]_{M} = \mathcal{L}_{\rho (u)} \pi \ ,
$$
while
$$
\frac{1}{2} \{\{\{ S_{\mathfrak{g}*} ,\rho\} ,\rho \} , u \} =
(\wedge^{2} \rho) (\gamma (u)) \ ,  
$$
and 
$$
\frac{1}{2} \{ \{\{\{ \Psi_{M}
,\rho \} , \pi\} ,\pi \} , u \}= (\wedge^2\pi^\sharp)(i_{\rho(u)} \Psi_M)  \ .
$$

Condition ($D$) is the relation
\begin{equation} \label{wedge3}
\frac{1}{2} [\pi ,\pi ]_{M} = (\wedge^{3}\rho) (\Psi_{\mathfrak{g}}) +
(\wedge^{3} \pi^{\sharp}) (\Psi_{M}) \ . 
\end{equation}

\subsection{Particular cases}\label{6.4} 
In the light of  
Proposition \ref{pi+rho} and formulas (\ref{rep}), (\ref{infpoisson}) and 
(\ref{wedge3}), we can interpret several important
particular cases of Poisson functions of the type $\pi + \rho$. 

\smallskip

\noindent
$\bullet$
{\bf Case $\rho = 0$}, already studied in section \ref{maurer}.
Conditions ($A$), ($B$) and ($C$) are identically satisfied and 
($D$) is the condition that $M$ be a {\it {twisted Poisson manifold}}.
If $\rho =0$ and $\Psi_{M}=0$, then ($D$) is the condition that
$M$ be a {\it{Poisson manifold}}.

\smallskip

\noindent
$\bullet$
{\bf Case $\Psi_{M} = 0$.} 
While condition ($A$) is identically satisfied, 
conditions ($B$), ($C$) and ($D$) 
express the fact that $M$ is a {\it quasi-Poisson $\mathfrak{g}$-space},
the version of the quasi-Poisson $G$-spaces
in the sense of \cite{AK}
in which only an infinitesimal Lie algebra action is assumed. When the
Lie group $G$ is connected and simply connected, conditions ($B$),
($C$) and ($D$) imply that $M$ is a {\it
  {quasi-Poisson $G$-space}}, and conversely.

\noindent
$\bullet$
{\bf Case $\Psi_M=0$ and $S_{\mathfrak{g}*} =0$.} Conditions ($B$), ($C$)
and ($D$) are\\
($B$) $M$ is a $\mathfrak{g}$-manifold,\\
($C$) $\pi$ is a $\mathfrak{g}$-invariant bivector,\\
($D$) $\frac{1}{2} [\pi ,\pi]_{M} = (\wedge^{3}\rho)
(\Psi_{\mathfrak{g}})$.
 
If $\Psi_{\mathfrak{g}}$ is the Cartan $3$-vector of the Lie
algebra $\mathfrak{g}$ of a connected and simply connected Lie group
with a bi-invariant scalar product, conditions ($B$), ($C$) and ($D$) 
express the fact that $M$ is a {\it quasi-Poisson $\mathfrak
  g$-manifold}, the version of the quasi-Poisson $G$-manifolds
in the sense of \cite{AKM} 
in which only an infinitesimal Lie algebra action is assumed.
When the
Lie group $G$ is connected and simply connected, conditions ($B$),
($C$) and ($D$) imply that $M$ is a {\it
  {quasi-Poisson $G$-manifold}}, and conversely.

\noindent
$\bullet$
{\bf Case $\Psi_{M}= 0$ and $\Psi_{\mathfrak{g}} = 0$.} In this
case, $(\mathfrak{g}, \mathfrak{g}^*)$ is a Lie bialgebra. 
Condition ($D$) expresses the fact that $\pi$ is a Poisson bivector, 
and equations (\ref{rep}) and (\ref{infpoisson}) show that conditions
($B$) and 
($C$) express the fact that $\rho$ is an 
{\it {infinitesimal Poisson action of the Lie bialgebra $(\mathfrak{g},
\mathfrak{g}^{*})$ on the Poisson manifold}} $M$
in the sense of Lu and Weinstein \cite{LW} \cite{L} (which can also be
called a {\it{Lie\break bialgebra action}}), corresponding to a
Poisson action of the connected and simply connected Poisson-Lie group
with Lie algebra $\mathfrak g$.

\medskip 

\noindent{\bf{Remark}}
The method described 
here for the characterization of Poisson and quasi-Poisson stuctures can
be used to recover conditions defining Poisson-Nijenhuis \cite{KM} and
Poisson-quasi-Nijenhuis \cite{SX} structures.

\subsection{The Lie algebroid 
structure of $V^*= T^*M \times \mathfrak g$}
Whenever $\sigma$ 
is a Poisson function with respect to a structure $S$ on $V$,
with $e^{-\sigma} S$,
$(V,V^{*})$ becomes a quasi-Lie bialgebroid.
Therefore when  
$\sigma = \pi + \rho$ 
is a Poisson function with respect to the structure 
$S = S_{\mathfrak{g}} + (S_{\mathfrak{g}^{*}} + S_{M})+
(\Psi_{\mathfrak{g}} + \Psi_{M})$ on $V= TM \times {\mathfrak g}^*$,
there is a Lie algebroid structure on $V^* =
T^{*}M
  \times \mathfrak{g}$, with anchor $\pi^{\sharp} + \rho$ and
Lie bracket
\begin{equation}\label{gamma}
\gamma_{\sigma} = S_{\mathfrak{g}} - \{ S_{\mathfrak{g}*} + S_{M} , \pi +\rho\}
+ \frac{1}{2} \{ \{ \Psi_{\mathfrak{g}} + \Psi_{M} ,\pi + \rho \},\pi
+\rho\} \ ,
\end{equation}
and $\{\gamma_{\sigma} , . \}$ is a differential on $\varGamma
(\wedge^{\bullet} (T M \times \mathfrak{g}^{*}))$.
Dually, there is a bracket $\mu_{\sigma}$ on $TM\times
\mathfrak{g}^{*}$, but the
Jacobi identity is not satisfied in general and the derivation
$\{\mu_\sigma, . \}$ on $\varGamma
(\wedge^{\bullet} (T^* M \times \mathfrak{g}))$
does not square to zero in general, since
$(V,V^{*})$ is only a quasi-Lie bialgebroid. From formula
(\ref{gamma}) and Proposition \ref{pi+rho}, we obtain: 
\begin{theorem}\label{liealgd}
When conditions ($A$)--($D$) are satisfied, 
$T^{*}M \times \mathfrak{g}$ is a Lie algebroid with anchor $\pi^{\sharp} + \rho$ and
Lie bracket  
\begin{equation}\label{bracket}
\begin{array}{llll}
\vspace{.1cm}
\gamma_\sigma & = & S_{\mathfrak{g}} - \{ S_{{\mathfrak{g}}^*},\rho\} - \{ S_{M} , \pi\} - \{
S_{M}, \rho\} \\
\vspace{.1cm}
& + & \frac{1}{2} \{\{ \Psi_{\mathfrak{g}},\rho\} ,\rho\} + 
\frac{1}{2} \{\{ \Psi_{M} ,\pi \} ,\pi\} 
+ \{\{ \Psi_{M} ,\pi \} ,\rho\} + \frac{1}{2} \{\{ \Psi_{M} ,\rho \} ,\rho\} \, .
\end{array}
\end{equation}
\end{theorem}
 
We shall now show that the preceding general formula
yields the brackets of \cite{L}, \cite{BC} and
\cite{BCS} as particular cases.

\smallskip

\noindent
{\bf Case $\rho = 0$.} Formula (\ref{bracket}) reduces to 
$ \gamma_\sigma 
= S_{\mathfrak{g}} - \{ S_{M} ,\pi \} + \frac{1}{2} \{ \{ \Psi_{M} ,\pi
\} ,\pi\}$. The Lie algebroid structure of $V^* = T^*M \times
\mathfrak g$ is the direct sum of the point-wise Lie bracket of
sections of $M \times \mathfrak g \to M$ and the Lie algebroid bracket
of \v Severa and Weinstein \cite{SW} on
$\varGamma (T^{*}M)$ for the twisted Poisson manifold $(M, \pi, \Psi_M)$.

\smallskip

\noindent
{\bf Case $\Psi_M = 0$.} Formula (\ref{bracket}) reduces to 
$$
\gamma_\sigma = S_{\mathfrak{g}} - \{ S_{\mathfrak{g}*} ,\rho \} - \{
S_{M},\pi \} -
\{ S_{M} ,\rho\}
+ \frac{1}{2} \{ \{ \Psi_{\mathfrak{g}}, \rho \} , \rho \}\ .
$$
For $u$, $v \in \varGamma (M \times \mathfrak{g})$ and 
$\alpha$, $\beta \in \varGamma (T^{*} M)$, we obtain the following
expressions entering in the brackets of sections of $T^*M \times \mathfrak
g$.
$$
\left\lbrace
\begin{array}{ll}
\vspace{0.2cm}
 \{\{u, S_{\mathfrak{g}} - \{ S_{M} ,\rho\}\},v\} =
[u,v]_{\mathfrak{g}} +
\mathcal{L}_{\rho(u)} v - \mathcal{L}_{\rho (v)} u \ ,\\
\vspace{0.2cm}
\{\{ \alpha, \{S_{\mathfrak{g}*},\rho \} \},u\} 
 = - i_{\rho^{*}(\alpha)} \{
S_{\mathfrak{g}^{*}},u \} ={\rm{ad}}^{*}_{\rho^{*}(\alpha)} u \ ,\\{}
\vspace{0.2cm}
\{\{ \alpha, \{S_{M} ,\pi \}\} ,u \} = -  \mathcal{L}_{\pi^\sharp(\alpha)} u
\ ,\\{}
\vspace{0.2cm}
\{\{ \alpha, \{S_{M} , \rho \}\} ,u \} = \mathcal{L}_{\rho(u)} \alpha
\ ,\\{}
\vspace{0.2cm}
\{\{ \alpha, \{S_{M} ,\pi \} \},\beta \} = - [\alpha ,\beta]_{\pi} 
\ ,\\{}
\frac{1}{2} \{\{\alpha, \{\Psi_{\mathfrak{g}},\rho \} ,\rho\}
\},\beta\} 
= i_{(\wedge^{2} \rho^{*})(\alpha \wedge \beta)}
\Psi_{\mathfrak{g}} \ ,
\end{array}
\right.
$$
where $\mathcal L$ denotes the Lie derivation of vector-valued
functions and of forms by vectors, and
${\rm{ad}}^*$ is defined by means of the bracket of ${\mathfrak g}^*$.
The bracket defined by $\gamma_\sigma$ is therefore 
$$
\left\lbrace
\begin{array}{ll}
\vspace{0.2cm}
[u,v] =
[u,v]_{\mathfrak{g}} +
\mathcal{L}_{\rho(u)} v - \mathcal{L}_{\rho (v)} u \ ,\\
\vspace{0.2cm}
[\alpha, u] =  \mathcal{L}_{\pi^\sharp(\alpha)} u
- \mathcal{L}_{\rho(u)} \alpha
- {\rm{ad}}^{*}_{\rho^{*}(\alpha)} u 
 \ ,\\{}
\vspace{0.2cm}
[\alpha,\beta ] = [\alpha ,\beta]_{\pi} + 
i_{(\wedge^{2} \rho^{*})(\alpha \wedge \beta)}
\Psi_{\mathfrak{g}} \ . 
\end{array}
\right.
$$
The bracket $[u,v]$ is the {\it {transformation Lie
  algebroid bracket}} \cite{M1} \cite{M2} on\break 
$M \times \mathfrak{g} \to M$.
Summarizing this discussion, we obtain
\begin{proposition}
If $\Psi_{M} = 0$, then $M$ is a quasi-Poisson $\mathfrak g$-space in
the sense of
{\rm {\cite{AK}}}
and
the Lie algebroid bracket of $T^{*} M\times
\mathfrak{g}$ is the bracket of Bursztyn, Crainic and \v Severa
{\rm{\cite{BCS}}}.
In particular, if $\Psi_{M} = 0$ and 
$S_{\mathfrak{g}^{*}} = 0$, then 
$M$ is a quasi-Poisson $\mathfrak g$-manifold in the sense of {\rm{\cite{AKM}}},
and
the Lie algebroid bracket of $T^{*} M\times
\mathfrak{g}$ is the bracket of Bursztyn and Crainic
{\rm {\cite{BC}}}.
\end{proposition}

\smallskip

\noindent
{\bf Case $\Psi_M = \Psi_{\mathfrak g}= 0$.} Formula (\ref{bracket}) reduces to 
$$
\gamma_\sigma =
S_{\mathfrak{g}} - \{ S_{\mathfrak{g}^{*}} , \rho \} - \{ S_{M} ,\pi\}
- \{ S_{M},\rho\} \ .
$$
Introducing the notations of Lu \cite{L}, 
the bracket of Bursztyn, Crainic and {\v S}evera reduces to the following 
expressions,
for $\alpha, \, \beta \in \varGamma (T^{*} M)$, and constant sections 
 $u, \, v$ of $M \times \mathfrak g$,
$$
\left\lbrace
\begin{array}{ll}
\vspace{.2cm}
[u , v] & = [u,v]_{\mathfrak{g}}\\{}
\vspace{.2cm}
[\alpha , u] & 
= D_{\alpha} u - D_{u} \alpha\\{}
\vspace{.2cm}
[\alpha ,\beta] & = [\alpha ,\beta ]_{\pi}
\end{array}
\right.
$$
\begin{proposition}
If $\Psi_{M} = 0$ and $\Psi_{\mathfrak g}=0$, then $M$ is a manifold
with a Lie bialgebra action and the Lie algebroid bracket of
$T^{*} M \times \mathfrak{g}$ is the bracket of Lu {\rm{\cite{L}}},
defining a matched pair of Lie  algebroids.
\end{proposition}

\subsection{The twisted differential}

Let us determine the differential $d_{\gamma_\sigma} = \{\gamma_\sigma, . \}$ on
$\varGamma(\wedge^\bullet (TM \times {\mathfrak g}^*))$, where
$\gamma_\sigma$ is defined by (\ref{bracket}).
The particular case of the quasi-Poisson $\mathfrak g$-spaces was
recently treated in \cite{BC2007}.

We first prove that the image of a section $X \otimes \eta$ of 
$\wedge^k TM \otimes \wedge^\ell {\mathfrak g}^*$ 
is a section of $\sum_{-1 \leq j \leq 2} 
\wedge^{k+j} TM \otimes \wedge^{\ell-j+1} {\mathfrak g}^*$.
We shall write $\varGamma({\mathfrak g}^*)$ for  $\varGamma(M \times 
{\mathfrak g}^* \to M)$. In fact, for $X \in \varGamma(\wedge^{k}
TM)$,
$$
\left\lbrace
\begin{array}{ll}
\vspace{.1cm}
&\{\{S_M,\pi\}, X\} \, \,
{\mathrm{and}}
\, \, \{\{\{ \Psi_{M} ,\pi \} ,\pi\}, X\} \in \varGamma(\wedge^{k+1} TM) \,
, \\
\vspace{.1cm}
&
\{\{S_M,\rho \}, X\} 
\, \,
{\mathrm{and}} \, \, 
\{ \{\{ \Psi_{M} ,\pi \} ,\rho\}, X\} \in 
\varGamma(\wedge^{k} TM
\otimes 
{\mathfrak g}^*) \, ,  \\
\vspace{.1cm}
&\{\{\{ \Psi_{M} ,\rho \} ,\rho\} , X \} \in \varGamma(\wedge^{k-1} TM
\otimes \wedge^2 
{\mathfrak g}^*) \, , 
\end{array}
\right.
$$
and for $\eta \in \varGamma(\wedge^{\ell} {\mathfrak g}^*)$,
$$
\left\lbrace
\begin{array}{ll}
\vspace{.1cm}
&
\{S_{\mathfrak g}, \eta\} \, \,
{\mathrm{and}} \, \, 
\{\{S_M,\rho\}, \eta\} 
 \in  
\varGamma(\wedge^{\ell+1} {\mathfrak g}^*)\, , \\
\vspace{.1cm}
& 
\{\{S_{{\mathfrak g}^*}, \rho \}, \eta\} \in
\, \,
{\mathrm{and}} \, \, \{\{S_M,\pi\}, \eta\} 
\in \varGamma(TM \otimes \wedge^{\ell} {\mathfrak g}^*) \ , \\
\vspace{.1cm}
&
\{ \{\{ \Psi_{\mathfrak{g}},\rho\} ,\rho\}  , \eta\} \in 
\varGamma(\wedge^{2} TM \otimes \wedge^{\ell - 1} {\mathfrak g}^*) \ ,
\end{array}
\right.
$$
while all other brackets vanish.

Each derivation is determined by
its values on the elements of degree $0$ and $1$. If $f \in C^{\infty}(M)$,
\begin{equation}\label{functions}
(d_{\gamma_\sigma}f) (\alpha + u) = 
(\pi^\sharp(\alpha) + \rho(u))\cdot f \ ,
\end{equation}
for $\alpha \in \varGamma(T^*M)$ 
and $u \in \varGamma({\mathfrak g})$.  
If $X \in \varGamma (TM)$, $d_{\gamma_\sigma} (X)$ is the sum of the
following terms,
$$
\left\lbrace
\begin{array}{ll}
\vspace{.1cm}
& - \{\{S_M,\pi\}, X\} + \frac{1}{2}\{\{\{ \Psi_{M} ,\pi \} ,\pi\}, X\} =
 [\pi,X]_M +(\wedge^2\pi^\sharp)(i_X\Psi_M) \\
\vspace{.1cm}
&  \quad \quad \quad \quad \quad \quad  \quad \quad \quad \quad \quad
\quad \quad \quad  \quad \quad \quad
\in \varGamma(\wedge^{2} TM) \, , \\
\vspace{.1cm}
& - \{\{S_M,\rho \}, X\}  +  \{\{\{ \Psi_{M} ,\pi \},\rho\}, X\} 
= [\rho(.), X]_M +
 (\pi^\sharp \wedge \rho)(i_X\Psi_M)
\\
\vspace{.1cm}
& \quad \quad \quad \quad \quad \quad  \quad \quad \quad \quad \quad
\quad \quad \quad  \quad \quad \quad
\in \varGamma(TM \otimes
{\mathfrak g}^*) \ , \\
\vspace{.1cm}
&
\frac{1}{2}\{\{\{ \Psi_{M} ,\rho \} ,\rho\}, X\} =
(\wedge^2\rho)(i_X\Psi_M) \in \varGamma( \wedge^2 
{\mathfrak g}^*) \ , 
\end{array}
\right.
$$
where $[\rho(.), X]_M : u \in {\mathfrak g} \mapsto [\rho(u), X]_M
\in \varGamma(TM)$. For $\eta \in \varGamma ({\mathfrak g}^*)$, 
$d_{\gamma_\sigma} (\eta)$ is the sum of the following terms,
\begin{equation}\label{eta}
\left\lbrace
\begin{array}{ll}
\vspace{.1cm}
& \{S_{\mathfrak g}, \eta\} - \{\{S_M,\rho\}, \eta\}
= d_\mathfrak g \eta + 
\ll {\mathcal L}_{\rho(.)} \eta, . \gg \, \in 
\varGamma(\wedge^2 {\mathfrak g}^*) \ , \\
\vspace{.1cm}
&
- \{\{S_{{\mathfrak g}^*}, \rho \}, \eta\} \! - \! 
\{\{S_M,\pi\}, \eta\}
= \rho({\rm{ad}}^*_\eta(.)) + {\mathcal
  L}_{\pi^\sharp (.)} \eta \in \varGamma(TM \otimes {\mathfrak g}^*), \\
\vspace{.1cm}
&
\frac{1}{2}
\{ \{\{ \Psi_{\mathfrak{g}},\rho\} ,\rho\}  , \eta\} =
- (\wedge^2\rho)(i_\eta \Psi_{\mathfrak g}) \in 
\varGamma(\wedge^2 TM) \ , 
\end{array}
\right.
\end{equation}
where 
$
\ll {\mathcal L}_{\rho(.)} \eta, . \gg : (u, v) \in \wedge^2
{\mathfrak g} \mapsto
\langle {\mathcal L}_{\rho(u)} \eta, v \rangle -
\langle {\mathcal L}_{\rho(v)} \eta, u \rangle \in C^\infty(M)$, 
$\rho({\rm{ad}}^*_\eta(.)) : u \in {\mathfrak g} \to
\rho({\rm{ad}}^*_\eta (u))\in \varGamma(TM)$,
and ${\mathcal
  L}_{\pi^\sharp (.)} \eta: \alpha \in \varGamma(T^*M) \mapsto {\mathcal
  L}_{\pi^\sharp (\alpha)} \eta \in \varGamma({\mathfrak g}^*)$.
The derivation $d_{\gamma_\sigma}$ is then extended to all 
sections of 
$\wedge^\bullet (TM \times {\mathfrak g}^*)$
by the graded Leibniz rule. We have thus obtained the following 
\begin{theorem} Let 
$\sigma = \pi + \rho$ be a Poisson function with respect to the structure
$S = S_{\mathfrak{g}} + (S_{\mathfrak{g}^{*}} + S_{M})+
(\Psi_{\mathfrak{g}} + \Psi_{M})$. 

\noindent(i) For $\gamma_\sigma$ defined by
(\ref{bracket}),
$d_{\gamma_\sigma}= \{\gamma_\sigma, . \}$ is a differential on 
$\varGamma(\wedge^\bullet (TM \times {\mathfrak g}^*))$.

\noindent(ii) $d_{\gamma_\sigma}= 
\sum_{-1 \leq j \leq2} d_{(j, 1-j)}$, where 
$$d_{(j,1-j)} : 
\varGamma(\wedge^k TM \otimes \wedge^\ell {\mathfrak g}^*) \to 
\varGamma(
\wedge^{k+j} TM \otimes \wedge^{\ell+1-j} {\mathfrak g}^*) \, ,
$$
and
$$
\begin{array}{ll} 
d_{(-1,2)} = \frac{1}{2}\{\{\{ \Psi_{M} ,\rho \} ,\rho\} , . \} \, ,\\
d_{(0,1)} = \{ - \{S_M,\rho \} +  \{\{ \Psi_{M} ,\pi \},\rho\} +
S_{\mathfrak g}, . \} \, ,\\
d_{(1,0)} =\{ - \{S_M,\pi\} + \frac{1}{2} \{\{ \Psi_{M} ,\pi \} ,\pi\} - 
 \{S_{{\mathfrak g}^*}, \rho \},  .\} \, ,\\  
d_{(2,-1)} = \frac{1}{2}
\{ \{\{ \Psi_{\mathfrak{g}},\rho\} ,\rho\}  , . \} \, .
\end{array}
$$
\noindent(iii)
For $f \in C^\infty(M)$ and $\eta \in \varGamma(M \times {\mathfrak
  g}^* \to M)$,
$d_{\gamma_\sigma}(f)$ and $d_{\gamma_\sigma} (\eta)$ are 
determined by Equations
(\ref{functions}) and (\ref{eta}) while, for $X \in \varGamma(TM)$, 
$$
d_{\gamma_\sigma} (X) =  [\pi , X]_M + [\rho(.), X]_M + 
( \wedge^2\pi^\sharp + \pi^\sharp
  \wedge \rho + \wedge^2 \rho ) (i_X\Psi_M) \ .
$$
\end{theorem}

These formulas simplify in each of the particular 
cases listed in Section \ref{6.4}. 
In the case of the quasi-Poisson $\mathfrak g$-spaces,
$d_{\gamma_\sigma} (X) =  [\pi,X]_M + {\mathcal L}_{\rho(.)} X$. 
From this formula and from (\ref{functions}), it follows that 
the restriction of $d_{\gamma_\sigma}$ to the space of $\mathfrak g$-invariant
multivectors on $M$ is the differential
of the {\it {quasi-Poisson cohomology}} 
introduced in \cite{AKM}. This fact was observed in \cite{T}. 

\medskip

\noindent{\bf Remark}
Thoughout this Section, the tangent bundle $TM$ can be
replaced by an arbitrary Lie algebroid over $M$, provided that
the de Rham differential is replaced by the differential associated
with the Lie algebroid in order to yield more general results.

\subsection*{Acknowledgments}
The main results of this paper were presented at the international
conference ``Higher structures in Geometry and Physics'' which was held
in honor of Murray Gerstenhaber and Jim Stasheff 
at the Institut Henri Poincar\'e in Paris in January 2007.
I am very grateful to the organizers, 
Alberto Cattaneo and Ping Xu, for the invitation to participate in
this exciting conference.
 
I thank Jim Stasheff, Murray Gerstenhaber, and Dmitry Roytenberg 
for their remarks and useful conversations.

%%%%%%%%%%%%%%%%%%%%%%%%%%%%%%%%%%%%%%%%%%%%%%%%%%%%%%%%%%%%%%%%%%%%%%%%%%

%%%%%%%%%%%%%%%%%%%%%%%%%

\begin{thebibliography}{AFMO}
%%%%%%%%%%%%%%%%%%%%
  

\bibitem{AK}
Alekseev A., Kosmann-Schwarzbach Y.:
Manin pairs and moment maps. J. Diff. Geometry, \textbf{56},
133--165 (2000)

\bibitem{AKM}
Alekseev A., Kosmann-Schwarzbach Y., Meinrenken E.:
Quasi-Poisson manifolds. Canadian J. Math., \textbf{54}, 3--29
(2002)


\bibitem{AX}
Alekseev A., Xu P.:
Derived brackets and Courant algebroids. Unpublished manuscript (2000)

\bibitem{B1}
Bangoura M.:
Alg\`ebres quasi-Gerstenhaber diff\'erentielles. Travaux math\'ematiques
(Luxembourg), \textbf{16}, 
299--314 (2005)

\bibitem{B2}
Bangoura M.:
Alg\`ebres d'homotopie associ\'ees \`a une proto-big\`ebre de Lie.
Canadian J. Math., \textbf{59}, 696--711 (2007)

\bibitem{BK}
Bangoura M., Kosmann-Schwarzbach Y.: 
The double of a Jacobian quasi-bialgebra.  
Lett. Math. Phys., \textbf{28}, 13--29 (1993)

\bibitem{BC}
Bursztyn H., Crainic M.:
Dirac structures, momentum maps, and
quasi-Poisson manifolds. In:
J. Marsden, T. Ratiu (eds),
The Breadth of Symplectic and Poisson
Geometry, Progr. Math., \textbf{232}, 1--40. Birkh\"auser, Boston (2005)

\bibitem{BC2007}
Bursztyn H., Crainic M.:
Dirac geometry, 
quasi-Poisson actions and $D/G$-valued moment maps, arXiv:0710.0639

\bibitem{BCS}
Bursztyn H., Crainic M., \v Severa P.:
Quasi-Poisson structures
as Dirac structures. Travaux Math\'ematiques (Luxembourg),
\textbf{16}, 41--52 (2005)

\bibitem{CDW}
Coste A., Dazord P., Weinstein A.:
Groupo\"{\i}des symplectiques. 
Publ. D\'ep. Math. Univ. Claude Bernard Lyon, Nouvelle S\'er. \textbf{2}/A,
1--62 (1987)

\bibitem{C}
Courant T.:
Dirac manifolds.  Trans. Amer. Math. Soc.,  \textbf{319}, 
631--661 (1990)

\bibitem{CW}
Courant T., Weinstein A.:
Beyond Poisson structures. 
Actions hamilto\-niennes de groupes. Troisi\`eme th\'eor\`eme de Lie,
S\'emin. Sud-Rhodan. G\'eom. VIII (Lyon, 1986), Trav. Cours, \textbf{27}, 39--49
(1988)

\bibitem{D}
Drinfeld V.:
Hamiltonian structures on Lie groups, Lie bialgebras and the geometric
meaning of classical Yang-Baxter equations.
Dokl. Akad. Nauk SSSR,  \textbf{268}, 285--287  (1983);
translation in Soviet Math. Doklady, \textbf{27}, 68--71 (1983)

\bibitem{D2}
Drinfeld V.:
Quasi-Hopf algebras. Algebra i Analiz, \textbf{1}, 114--148 (1989); 
translation in Leningrad Math. J., \textbf{1}, 1419--1457 (1990)

\bibitem{E1}
Ehresmann C.:
Cat\'egories topologiques et cat\'egories diff\'erentiables.
Centre Belge Rech. Math., Colloque G\'eom. Diff\'er. Globale
(Bruxelles 1958), 137--150 (1959) 

\bibitem{E2}
Ehresmann C.:
Sur les cat\'egories diff\'erentiables. 
Atti Convegno internaz. Geom. diff. (Bologna 1967),
31--40 (1970) 

 
\bibitem{G}
Gerstenhaber M.: 
The cohomology structure of an associative ring.  
Ann. of Math. (2),  \textbf{78}, 267--288  (1963) 

\bibitem{GG}
Gerstenhaber M., Giaquinto A.:
Boundary
solutions of the classical Yang--Baxter equation. Lett. Math. Phys.,
{\textbf{40}}, 337--353 (1997)
 
\bibitem{HY}
Hodges T. J., Yakimov M.: 
Triangular Poisson structures on Lie groups and symplectic
  reduction. In: 
Noncommutative Geometry and Representation Theory in Mathematical
Physics,  
Contemp. Math., \textbf{391}, 123--134. Amer. Math. Soc., 
Providence, R.I. (2005) 

\bibitem{H1}
Huebschmann J.:
Poisson cohomology and quantization.  J. Reine Angew. Math.,
\textbf{408},
57--113 (1990) 

\bibitem{H2}
Huebschmann J.:
Higher homotopies and Maurer-Cartan algebras: quasi-Lie-Rinehart,
Gerstenhaber, and Batalin-Vilkovisky algebras. In: 
J. Marsden, T. Ratiu (eds), 
The Breadth of
Symplectic and Poisson Geometry, Progr. Math., \textbf{232}, 237--302. 
Birkh\"auser, Boston (2005)

\bibitem{J} 
Jacobson, N.:
On pseudo-linear transformations.
Proc. Natl. Acad. Sci. USA, \textbf{21}, 667--670 (1935)

\bibitem{KS} 
Klim\v c\'{\i}k C., Strobl T.:
WZW-Poisson manifolds.  J. Geom. Phys.,  \textbf{43}, 341--344  (2002)

\bibitem{yks92} 
Kosmann-Schwarzbach Y.:
Jacobian quasi-bialgebras and quasi-Poisson Lie groups. 
In: Mathematical Aspects of Classical Field Theory (Seattle 1991),  
Contemp. Math., \textbf{132},  459--489. Amer. Math. Soc., 
Providence, R.I. (1992) 

\bibitem{yks95} 
Kosmann-Schwarzbach Y.:
Exact Gerstenhaber 
algebras and Lie bialgebroids. Acta Appl. Math., \textbf{41}, 153--165 (1995)

\bibitem{yks96} 
Kosmann-Schwarzbach Y.:
From Poisson algebras to Gerstenhaber algebras.
Ann. Inst. Fourier (Grenoble), \textbf{46}, 1243--1274 (1996) 

\bibitem{yks04} 
Kosmann-Schwarzbach Y.:
Derived brackets. Lett. Math. Phys., \textbf{69}, 61--87 (2004)

\bibitem{yks05} 
Kosmann-Schwarzbach Y.:
Quasi, twisted, and 
all that$\ldots$ in Poisson 
geometry and Lie algebroid theory. In: J. Marsden, T. Ratiu (eds), 
The Breadth of Symplectic and
Poisson Geometry, Progr. Math., \textbf{232}, 363--389. Birkh\"auser, Boston
(2005) 

\bibitem{KL} 
Kosmann-Schwarzbach Y., Laurent-Gengoux C.:
The modular class of a twisted
Poisson structure. Travaux Math\'ematiques (Luxembourg), \textbf{16},
315--339 (2005) 


\bibitem{KMc}
Kosmann-Schwarzbach Y., Mackenzie K. C. H.:
Differential operators and actions of Lie algebroids. In: 
T. Voronov (ed), Quantization,
Poisson Brackets and Beyond, Contemp. Math., \textbf{315},
213--233. Amer. Math. Soc., Providence, R.I. (2002)

\bibitem{KM}
Kosmann-Schwarzbach Y., Magri F.:
Poisson-Nijenhuis structures. Ann. Inst. Henri Poincar\'e, S\'erie A,
\textbf{53},
35--81 (1990)

 
\bibitem{KY} 
Kosmann-Schwarzbach Y., Yakimov M.:
Modular classes of regular twisted Poisson structures on Lie
algebroids. Lett. Math. Phys., \textbf{80}, 183--197 (2007) 

\bibitem{KSt}
Kostant B., Sternberg S.:
Symplectic reduction, BRS cohomology, and infinite-dimensional
Clifford algebras.
Ann. Physics, \textbf{176}, 49--113 (1987)



\bibitem{Ko}
Koszul J.-L.: 
Crochet de Schouten--Nijenhuis et cohomologie. 
The Mathematical heritage of \'Elie Cartan (Lyon, 1984).
Ast\'erisque,  num\'ero hors s\'erie, 257--271 (1985)

\bibitem{LR}
Lecomte P., Roger C.:
Modules et cohomologies des big\`ebres de Lie.
C. R. Acad. Sci. Paris S\'er. I Math.,  \textbf{310}, 405--410 (1990)

\bibitem{Li}
Lichnerowicz A.:
Les vari\'et\'es de Poisson et leurs alg\`ebres de Lie associ\'ees.
J. Differential Geom., \textbf{12}, 253--300 (1977)

\bibitem{LWX}
Liu Z.-J., Weinstein A., Xu P.:
Manin triples for Lie bialgebroids. 
J. Differential Geom., \textbf{45},  547--574 (1997)

\bibitem{Lo}
Loday J.-L.:
Une version non commutative des alg\`ebres de Lie: les alg\`ebres de
Leibniz. Enseign. Math., \textbf{39}, 269--293 (1993)

\bibitem{L}
Lu J.-H.:
Poisson homogeneous spaces and 
Lie algebroids associated to Poisson actions. 
Duke Math. J., \textbf{86}, 261--304 (1997) 

\bibitem{LW}
Lu J.-H., Weinstein A.:
Poisson Lie groups, dressing transformations, and Bruhat
decompositions.  J. Differential Geom., 
\textbf{31}, 501--526  (1990) 

\bibitem{M1} 
Mackenzie K. C. H.:
Lie groupoids and Lie Algebroids in Differential Geometry. London
Mathematical Society Lecture Note Series, \textbf{124}. Cambridge University
Press, Cambridge (1987)

\bibitem{M0} 
Mackenzie K. C. H.:
Lie algebroids and Lie pseudoalgebras. 
Bull. Lond. Math. Soc., \textbf{27}, 97--147 (1995)

\bibitem{M2} 
Mackenzie K. C. H.:
General Theory of Lie Groupoids and Lie Algebroids. London
Mathematical Society Lecture Note Series, \textbf{213}. Cambridge University
Press, Cambridge (2005)

\bibitem{MX} 
Mackenzie K. C. H., Xu P.:
Lie bialgebroids and Poisson groupoids. Duke Math. J., \textbf{73},
415--452 (1994)

\bibitem{MM}
Moerdijk I., Mr\v cun J.:
Introduction to Foliations and Lie Groupoids. Cambridge Studies in
Advanced Mathematics, \textbf{91}. Cambridge University Press, Cambridge (2003)

\bibitem{P}
Pradines J.:
Th\'eorie de Lie pour les groupo\"{\i}des diff\'erentiables.
Calcul diff\'erentiel dans la cat\'egorie des groupo\"{\i}des
infinit\'esimaux.
C. R. Acad. Sci. Paris S\'er. A-B, \textbf{264}, A245--A248 (1967)


\bibitem{R} 
Roytenberg D.: 
Quasi-Lie bialgebroids and twisted Poisson
manifolds.
Lett. Math. Phys., \textbf{61}, 123--137 (2002) 

\bibitem{Roytenberg} 
Roytenberg D.: 
e.mail message (2007)

\bibitem{S}
Saksida P.:
Lattices of Neumann oscillators and Maxwell-Bloch equations.
Nonlinearity,  \textbf{19}, 747--768  (2006)

\bibitem{SW}
\v Severa P., Weinstein A.:
Poisson geometry 
with a 3-form background. In: Noncommutative Geometry 
and String Theory (Yokohama, 2001), Progr. Theoret. Phys. 
Suppl., \textbf{144}, 145--154 (2001) 


\bibitem{Stasheff}
Stasheff J.:
Constrained Hamiltonians, BRS and homological algebra. In: 
Proceedings of the 
Conference on Elliptic Curves and Modular Forms in Algebraic Topology
(Princeton, 1986), Springer Lecture Notes in Math.,
\textbf{1326}, 150--160. Springer, Berlin (1988)

\bibitem{Sta}
Stasheff J.:
Differential graded Lie algebras, quasi-Hopf algebras and higher
homotopy algebras. In: Quantum Groups (Leningrad, 1990),  
Lecture Notes in Math., \textbf{1510}, 120--137. Springer, Berlin (1992)

\bibitem{SX}
Sti\'enon M., Xu P.:
Poisson quasi-Nijenhuis manifolds.
Comm. Math. Phys., \textbf{270}, 709--725 (2007)

\bibitem{Sto} 
Stolin A.: 
On rational solutions of Yang--Baxter 
equation for $\sl(n)$. Math. Scand., {\textbf{69}}, 57--80 (1991)

\bibitem{T} 
Terashima Y.:
On Poisson functions. J. Sympl. Geom., \textbf{6}, no. 1, 1--7 (2008)

\bibitem{V} 
Vaintrob A.: 
Lie algebroids and homological vector fields.
Uspekhi Mat. Nauk, \textbf{52}, 2(314),
161--162 (1997); translation in Russ. 
Math. Surv., \textbf{52}, 428--429 (1997) 

\bibitem{Vo} 
Voronov T.:
Graded manifolds and Drinfeld doubles for Lie
bialgebroids. In: T. Voronov (ed.), 
Quantization, Poisson Brackets and
Beyond, Contemp. Math., \textbf{315}, 131--168.
Amer. Math. Soc., Providence, R.I. (2002)

\bibitem{X}
Xu P.:
Gerstenhaber algebras and BV-algebras in Poisson geometry. 
Comm. Math. Phys., \textbf{200}, 545--560 (1999)

%%%%%%%%%%%%%%%%%%%%%%%%%%%%%%%%%%%%%%%%%%%%%%%%%%%%%%
\end{thebibliography}
\end{document}